\documentclass[11pt]{amsart}
\usepackage{amsmath,amssymb,amsfonts}
\usepackage{graphicx}

\newtheorem{defi}{Definition}[section]
\newtheorem{definition}[defi]{Definition}
\newtheorem{proposition}[defi]{Proposition}
\newtheorem{theorem}[defi]{Theorem}
\newtheorem{remark}[defi]{Remark}
\newtheorem{corollary}[defi]{Corollary}
\newtheorem{lemma}[defi]{Lemma}

\newcommand{\bb}[2]{\genfrac{[}{]}{0pt}{}{#1}{#2}}
\newcommand{\ka}{\kappa}

\newcommand{\la}{\lambda} 
\newcommand{\e}{\epsilon} 
\newcommand{\lep}{\lambda_k^{\epsilon}} 
\newcommand{\ve}{v^{\epsilon}} 
\newcommand{\dk}{d_{2\kappa}}

\newcommand{\R}{{\mathbb{R}}} 
\newcommand{\C}{{\mathbb{C}}}

\newcommand{\Gmn}{{\mathcal{G}_{m,n}}}

\newcommand{\On}{{\operatorname{O}(n,\R )}}
\newcommand{\Om}{{\operatorname{O}(m,\R )}}
\newcommand{\Onm}{{\operatorname{O}(n-m,\R )}}

\newcommand{\Glm}{{\operatorname{GL}(m,\R)}}

\newcommand{\s}{\sigma}

\newcommand{\ki}{\kappa^{(i)}}
\newcommand{\kii}{\kappa_{(i)}}

\newcommand{\kjj}{\kappa_{(j)}}

\newcommand{\kkii}{{\kappa'}_{(i)}}
\newcommand{\kkj}{{\kappa'}^{(j)}}

\newcommand{\Pk}{P_{\kappa}}
\newcommand{\Qk}{Q_{\kappa}}

\newcommand{\Pmu}{P_{\mu}}
\newcommand{\Hka}{\operatorname{H}_{m,n}^{2\kappa}}
\newcommand{\Pnu}{P_{\nu}}
\newcommand{\Ck}{C_{\kappa}}

\newcommand{\Cmu}{C_{\mu}}
\newcommand{\Cnu}{C_{\nu}}

\newcommand{\PP}{{\mathbb{P}}}

\newcommand{\Harm}{\operatorname{Harm}}

\newcommand{\Id}{\operatorname{Id}}

\newcommand{\pr}{\operatorname{pr}} 
\newcommand{\one}{{\bf 1}}

\begin{document}
\title{Linear Programming bounds for codes in Grassmannian spaces}
\author{Christine Bachoc}

\thanks{C. Bachoc is with the Laboratoire A2X, Institut de Math\'ematiques de
  Bordeaux, 351 cours de la Lib\'eration, 33405 Talence, France, bachoc@math.u-bordeaux1.fr}

\begin{abstract} We develop the linear programming method in order 
to obtain bounds
for the cardinality of  Grassmannian codes endowed with the chordal distance.
We obtain a bound and its asymptotic version that generalize the well-known 
bound for codes in the real projective space obtained by Kabatyanskiy and 
Levenshtein, and improve the Hamming bound for sufficiently large
minimal distances.
\end{abstract}

\maketitle

\section{Introduction}

Philippe Delsarte has introduced the so-called
{\em linear programming method}, in order  to find bounds 
for the size of codes with prescribed
minimal distance, in the classical
case of codes over finite fields. 
This method, also called {\em Delsarte method}
or {\em polynomial method}, exploits a certain family
of orthogonal polynomials attached to the situation, the Krawtchouk 
polynomials, and their positivity property. These polynomials and
their properties are intimately related to the action of
the symmetric group on the Hamming space.
Delsarte method has proved to be very powerful,
and was extended to many other situations, where 
the underlying space is symmetric of rank one, and is
homogeneous under the action of a certain group of transformations.
Examples of such spaces are: the Johnson space, the Grassmannian space 
over a finite field,
the unit sphere of the Euclidean space, the projective spaces
over the real, complex and quaternionic fields (for these last spaces
see \cite{DGS}, \cite{DGS2}).

In recent years, codes over the real Grassmannian space
have attracted attention, motivated by their application to information theory,
more precisely to the so-called {\em space-time codes}, used for
multi-antenna
systems of communication. The distance usually
considered is the chordal distance, introduced in \cite{CHS}, and
defined in the following way (more details are given in the next subsection):
The Grassmannian space of $m$-dimensional subspaces 
of $\R^n$, where $m\leq n/2$, is denoted by $\Gmn$; to a pair $(p,q)$ of elements  of $\Gmn$ is associated
$m$ principal angles $\theta_1,\dots,\theta_m\in [0,\pi/2]$. Let
$y_i:=\cos^2\theta_i$. Then
\begin{equation*} 
d_c(p,q):=\sqrt{\sum_{i=1}^m \sin^2 \theta_i}=\sqrt{m-\sum_{i=1}^m y_i}.
\end{equation*}

In \cite{CHS}, the authors give bounds for the size of Grassmannian codes,
called the simplex and orthoplex bounds. The main drawback of these bounds 
is that they are only valid in a certain range of minimal distances.
In \cite{BN}, an asymptotic bound, derived from  the Hamming bound,
is given. Another approach is developed in 
\cite{BBC}, where bounds are given for codes whose principal angles
are subject to certain constraints (the so-called $f$-codes), which 
arise naturally from the notion of Grassmannian designs introduced
in \cite{BCN}.

In this paper, we extend Delsarte method to the Grassmannian codes, 
exploiting 
the zonal polynomials attached to $\Gmn$. These are symmetric polynomials 
in the $m$ variables $y_1,\dots, y_m$; they 
 belong to the family of orthogonal {\em generalized Jacobi
   polynomials} (see the next subsection).
In the second section, we recall, or settle the properties of these polynomials
needed to perform linear programming bounds; these properties are 
easy to obtain by straightforward generalization of the arguments used in the 
classical cases. In fact, the principles underlying the LP method would
remain true for the zonal polynomials attached to any symmetric space. The 
real difficulties start when one wants to actually perform explicit bounds,
because the polynomials have (for $m\geq 2$) several variables.
The low degree cases are still easy to manage; this is done in section \ref{sec3},
where we recover the simplex bound as the bound arising from the case of
degree one, and give new bounds from polynomials of degree $2$ and
$3$. 
In the forth section, we propose a strategy based on 
the eigenvalues of certain symmetric endomorphisms, which extends the one
variable method based on the zeros of the polynomials and on
Christoffel-Darboux formula, but avoids to deal with zeros of
polynomials in several variables. We obtain an upper bound for the
size of a code $C$ with minimal distance $\delta$, which is expressed
in terms of the largest eigenvalue (Theorem \ref{Teps} and Corollary \ref{ineg}).
Section \ref{sec5} settles the asymptotic behavior of this largest eigenvalue
(Theorem \ref{lambda}),
and in section \ref{sec6} we derive the following asymptotic version of the
bound:

\medskip
\begin{theorem}\label{Tass}
Let $C$ be a code in $\Gmn$ with minimal chordal distance $\delta$,
let $s:=m-\delta^2\in ]0,m[$ and let 
$$\rho:=\frac{m}{2}(-1+(1-\frac{s}{m})^{-1/2}).$$
Then, when $n\to +\infty$,

\begin{equation}\label{ass}
\frac{1}{n}\log |C|\lesssim m\big((1+\rho)\log(1+\rho)-\rho\log(\rho)).
\end{equation}
\end{theorem}
\medskip

Our bound coincides with the bound given by G. Kabatiansky and
V. Levenshtein in \cite{KL} for the case of the real projective space,
corresponding to $m=1$. But it beats the Hamming bound of \cite{BN}
only when the minimal distance is relatively big.

\subsection{Basic facts about Grassmannian spaces and their zonal
  polynomials.}

We repeat here, without proofs, some well-known facts about Grassmannian spaces and
their zonal polynomials. Some useful references for the mathematical background are: \cite{FH}, \cite{GW} for the
representations of the orthogonal group, \cite{JC}, \cite{GL},
 for the Grassmannian spaces and harmonic analysis on it, \cite{La1},
 \cite{La2}, \cite{Macdo}, \cite{D} for multivariate orthogonal polynomials.

The real Grassmannian space, denoted by
$\Gmn$ ($m\leq n/2$), is the set of $m$-dimensional $\R$-linear subspaces of $\R^n$. 
The orthogonal group $\On$ acts transitively on $\Gmn$; a transformation stabilizing
a given element $p_0$ also stabilizes its orthogonal complement
$p_0^{\perp}$ and therefore the stabilizer of $p_0$ is isomorphic to the direct product
$\Om\times \Onm$. Hence we derive the identification of $\Gmn$ with
the set of classes:
\begin{equation*}
\Gmn\simeq \On/\big(\Om\times \Onm\big)
\end{equation*}
from which $\Gmn$ inherits the structure of a (compact) dif\-fe\-ren\-tial
variety, and a $\On$-invariant measure that will be normalized so that
$\int_{\Gmn} dp=1$. It is worth noticing that the case $m=1$
corresponds to the real projective space.

In order to understand the action of $\On$ on pairs $(p,q)\in \Gmn^2$,
we need
to introduce the {\em principal angles} between $p$ and $q$.
These are $m$ angles $\theta_1,\dots, \theta_m\in [0,\pi/2]$ defined
in the following way:

Let $p_1\subset p$, $q_1\subset q$ be two lines such that the angle
$\theta_1$ between $p_1$ and $q_1$ is minimal. If $m=1$ we have
finished, otherwise
let $p'$ be the orthogonal complement of $p_1$ in $p$,
$q'$ be the orthogonal complement of $q_1$ in $q$; we define 
recursively $\theta_2,\dots,\theta_m$ to be the principal angles
associated
to the pair $(p',q')$ in ${\mathcal G}_{m-1,n}$.
We introduce the notation $y_i:=\cos^2\theta_i$; when needed, we may
denote rather $y_i(p,q)$, $\theta_i(p,q)$. A classical result on the
geometry of $\Gmn$ is the following:

\medskip
\begin{proposition}\label{orb}
Two pairs $(p,q)$ and $(p',q')$ are in the same orbit under the action
of the group $\On$, i.e. there exists $\sigma\in \On$ such that
$\sigma(p)=p'$ and $\sigma(q)=q'$, if and only if 
$$y_i(p,q)=y_i(p',q') \text{ for all } 1\leq i\leq m.$$
\end{proposition}
\medskip

The previous proposition expresses the fact that the orbits under the
action
of $\On$ of the pairs $(p,q)\in \Gmn^2$ are characterized by the
$m$-tuple of real numbers $(y_1(p,q),\dots,y_m(p,q))$. It becomes
clear
that, for $m\geq 2$, $\Gmn$ is not $2$-point homogeneous, i.e. a
single distance on $\Gmn$ cannot characterize these orbits
(while it is the case for other spaces of interest in coding theory,
like
the Hamming and binary Johnson spaces, or the unit sphere of the
Euclidean space). It is the reason why we shall deal with mutivariate
polynomials.
Also, it shows that the choice of a distance on $\Gmn$ is sort of 
arbitrary. We shall stick to the {\em chordal distance} in this paper,
as introduced in \cite{CHS}:
\begin{equation*}
d_c(p,q):=\sqrt{\sum_{i=1}^m \sin^2 \theta_i}=\sqrt{m-\sum_{i=1}^m y_i}.
\end{equation*}

Other possibilities are the Riemannian distance $\sqrt{\sum_{i=1}^m
  \theta_i^2}$
which behaves somewhat badly because it is not smooth; the max
  distance $\max_{i} \theta_i$, etc.. The ``product distance'' (which is not a distance
in the metric sense) $\big(\prod_i \sin\theta_i\big)$ seems to be relevant in the
  context
of space time codes. 

\smallskip
Now we consider the space $L^2(\Gmn)$ of functions $f:\Gmn\to \C$ such
that $\int_{\Gmn} |f|^2dp<+\infty$. This is a $\C$-vector space,
endowed
with the hermitian product:
$$<f,g>=\int_{\Gmn} f(p)\overline{g(p)}dp$$
and with the left  action of the orthogonal group given by:
$$(\sigma\cdot f)(p)=f(\sigma^{-1}(p))$$
(for which the above hermitian product is of course  invariant).

Its associated {\em algebra of zonal functions} (also called the Hecke
algebra) is:
\begin{align*}
{\mathcal Z}:=\{Z:&\Gmn^2\to \C \mid \ Z(p,\cdot),\ 
Z(\cdot,q)\in L^2(\Gmn)\text{ and }\\
&Z(\sigma(p),\sigma(q))=Z(p,q)\text{ for all }\sigma\in \On\}
\end{align*}
Form Proposition \ref{orb}, since $Z\in {\mathcal Z}$ is constant on
the orbits of $\On$ on $\Gmn^2$, it can be given the form:
$Z(p,q)=z(y_1(p,q),\dots,y_m(p,q))$ for some function $z$.

The explicit decomposition into $\On$-irreducible subspaces of\\
$L^2(\Gmn)$,
and the corresponding structure of $\mathcal Z$, where investigated
for the first time by James and Constantine (\cite{JC}). It is now a
standard
result on the representation of the classical groups (see \cite{GW}).

Recall that the irreducible representations of $\On$ are (up to a power of the
determinant)
naturally indexed by partitions
$\kappa=(\kappa_1,\dots,\kappa_n)$, where $\kappa_1\geq \dots\geq
\kappa_n\geq 0$ (we may omit the last parts if they are equal to $0$).
Following \cite{GW}, let them be denoted by $V_n^{\kappa}$.
For example, $V_n^{()}=\C\one$, and $V_n^{(k)}=\Harm_k$ the space 
of homogeneous of degree $k$,  harmonic polynomials in $n$ variables.

The length $\ell(\kappa)$ of a partition $\kappa$ is the
number of its non zero parts, and its degree 
$\deg(\kappa)$ also denoted by $|\kappa|$ equals $\sum_{i=1}^n\kappa_i$.

\smallskip
Then, the decomposition of $L^2(\Gmn)$ is as follows:

\begin{equation*}
L^2(\Gmn)\simeq \oplus V_n^{2\kappa}
\end{equation*}
where $\kappa$ runs over the partitions of length at most $m$ and
$2\kappa$ stands for $(2\kappa_1,\dots,2\kappa_m)$, meaning that
only partitions with even parts enter the decomposition.
We can see that the multiplicities in this decomposition are all equal
to one, which translates the fact that the space $\Gmn$ is a symmetric
space. Consequently, to each irreducible component $V_n^{2\kappa}$ is
associated a uniquely determined (up to a normalizing factor) 
zonal function $\Pk(y_1,\dots,y_m)$,
in the sense that 
\begin{equation*}
p \ (\text{resp. } q)\mapsto \Pk(y_1(p,q),\dots,y_m(p,q))\in V_n^{2\kappa}
\end{equation*}
and 
\begin{equation*}
{\mathcal Z}=\oplus_{\kappa, \ell(\kappa)\leq m}\C \Pk.
\end{equation*}

It turns out that the $\Pk$ are symmetric polynomials in the
$m$ variables $y_1,\dots,y_m$, of degree $|\kappa|$,
with rational coefficients once they are normalized by the condition  $\Pk(1,\dots,1)~=~1$.
Moreover, the set $(\Pk)_{|\kappa|\leq k}$ is a basis of the space 
of symmetric polynomials in the variables $y_1,\dots,y_m$ of degree at most
equal to $k$, denoted by $S_k$.

Since the irreducible subspaces of $L^2(\Gmn)$ are pairwise non
isomorphic, they are orthogonal for the $\On$-invariant hermitian
product defined above. This hermitian product
induces an hermitian product on the space of symmetric polynomials,
denoted by  $[,]$, for which the polynomials $\Pk$ are orthogonal.
More precisely, it is given by the positive measure, calculated in \cite{JC},

$$d\mu=
\lambda \prod
_{\substack{i,j=1\\i<j}}^m|y_i-y_j|\prod_{i=1}^my_i^{-1/2}(1-y_i)^{n/2-m-1/2} dy_i
$$
(where $\lambda$ is chosen so that $\int_{[0,1]^m}d\mu(y)=1$).
and 

\begin{equation*}
[f,g]=\int_{[0,1]^m}f(y)\overline{g(y)}d\mu(y).
\end{equation*}

One recognizes a special case of the orthogonal measure associated to {\em generalized
  Jacobi polynomials} (\cite{La2}).

We let  $\Pi_k$ be the subspace of $S_k$ generated by the polynomials
$\Pk$, with $|\kappa|=k$, so that we have the orthogonal decomposition:

\begin{equation*}
S_k=S_{k-1}\perp \Pi_k.
\end{equation*}
Let  the dimensions of $S_k$, $\Pi_k$ be denoted respectively by $s_k$, $\pi_k$.
The number $\pi_k$ is also equal to the number
of partitions $\kappa$ of $k$ in at most $m$ parts. These dimensions
 also depend on $m$, although it does not reflect on our notation, for the 
sake of simplicity.

In view of the explicit calculation of the polynomials $\Pk$, it is better
to use the following characterization, which involves the polynomials
$\Ck$, which are themselves the zonal polynomials associated to the symmetric space
$\Glm/\Om$ (these polynomials are Jack polynomials, normalized by $\Ck(1,\dots,1)=1$,
see \cite{JC}, \cite{Macdo}), and the differential operator $\Delta$
induced on $\C[y_1,\dots,y_m]^{S_m}$ by the Laplace Beltrami operator
of $\Gmn$. The
condition: for all $1\leq i\leq m$, $\kappa_i\geq \mu_i$ is denoted by: $\kappa\geq \mu$.

\begin{enumerate}
\item $\Pk$ is an eigenvector for the operator 

\begin{equation}\notag
\begin{split}
\Delta:=&\sum_{i=1}^m y_i^2\frac{\partial^2}{\partial y_i^2}
+\sum_{i\neq j=1}^m y_i^2(y_i-y_j)^{-1}\frac{\partial}{\partial y_i}\\
&+(\frac{n}{2}-m+1)\sum_{i=1}^m y_i\frac{\partial}{\partial y_i}-\sum_{i=1}^m y_i\frac{\partial^2}{\partial y_i^2}\\
&-\sum_{i\neq j=1}^m y_i(y_i-y_j)^{-1}\frac{\partial}{\partial y_i}
-\frac{1}{2}\sum_{i=1}^m \frac{\partial}{\partial y_i}
\end{split}
\end{equation}

\item $\Pk=\beta_{\kappa}\Ck+\sum_{\mu \mid \kappa>\mu }\beta_{\kappa,\mu}\Cmu$

\item $\Pk(1,\dots,1)=1$. 
\end{enumerate}
Condition (ii) is needed
to avoid the multiplicities of the  operator $\Delta$.

\medskip
\noindent{ \bf Examples:} the effective computation of the polynomials
$\Pk$ following the method described above leads to, up to the
normalization imposed by (iii):

\begin{align*}
&P_0 =1 \\
&P_{(1)} =s_1 -\frac{m^2}{n}\\
&P_{(11)} =\sigma_1 -\frac{(m-1)^2}{n-2}s_1+\frac{m^2(m-1)^2}{2(n-1)(n-2)} \\
&P_{(2)}=s_2 +\frac{2}{3}\sigma_1 -\frac{2(m+2)^2}{3(n+4)}s_1+\frac{m^2(m+2)^2}{3(n+2)(n+4)}
\end{align*}
where $s_1=\sum_{1\leq i\leq m} y_i$, $s_2=\sum_{1\leq i\leq m} y_i^2$,\\
$\sigma_1=\sum_{1\leq i<j\leq m} y_iy_j$.

\medskip
\noindent{\bf Remark:} The complex Grassmannian $\Gmn(\C)$ is more commonly used in the context of
space-time coding. It affords the transitive action of the unitary
group $U(\C,n)$; similarly one defines principal angles
$(\theta_1,\dots,\theta_m)$ between two elements of $\Gmn(\C)$. The
$U(\C,n)$ decomposition of $L^2(\Gmn(\C))$ and the associated zonal
polynomials are computed in \cite{JC} so one can play the same game
concerning bounds of codes. On the other hand, a bound is also obtained from the
embedding
$\Gmn(\C)\subset {\mathcal G}_{2m,2n}(\R)$ (if
$(\theta_1,\dots,\theta_m)$ are the principal angles
associated to a pair $(p,q)$ of elements in $\Gmn(\C)$, the $2m$
principal angles associated to the pair $(p,q)$, seen as elements of
${\mathcal G}_{2m,2n}(\R)$,
are simply $(\theta_1,\theta_1,\theta_2,\theta_2,\dots)$).

\section{Zonal polynomials associated to $\Gmn$ and the LP bound}\label{sec2}

In this section, we settle the properties of the polynomials $\Pk$
relevant for the LP bound, settle this bound, and show how the
Christoffel-Darboux
formula can be exploited in that context.

\smallskip
The dimension of $V_n^{\kappa}$ is
denoted by $d_{\kappa}$. Explicit formulas for $d_{\kappa}$ can be
found in \cite{FH}; however we do not need them before Section 5.

\medskip
\begin{proposition}\label{PP} The polynomials $\Pk$, normalized 
by the condition \\ $\Pk(1,\dots,1)=1$, satisfy:
\begin{enumerate}
\item $[\Pk,\Pk]=d_{2\kappa}^{-1}$
\item (Positivity property): For all finite set $C\subset\Gmn$,
\begin{equation*}
\sum_{p,q\in C}\Pk(y_1(p,q),\dots,y_m(p,q))\geq 0
\end{equation*}
\item Let $p_{\kappa,\mu}^{\nu}$ be defined by the property:
\begin{equation*}
\Pk\Pmu=\sum_{\nu}p_{\kappa,\mu}^{\nu}\Pnu
\end{equation*}
The numbers $p_{\kappa,\mu}^{\nu}$ are non-negative numbers.
\end{enumerate}
\end{proposition}
\medskip

\noindent\proof These properties where already pointed out in \cite{BBC}[Lemma 2.2]
and step on very general arguments (see \cite{K}[Theorem 3.1]). For the sake of completeness,
we briefly recall the arguments. Let $e_1,\dots,e_{d_{2\kappa}}$ be 
any orthonormal basis
of the subspace $\Hka$ of $L^2(\Gmn)$ isomorphic to $V_n^{2\kappa}$.
Let $\tilde{\Pk}(p,q):=\Pk(y_1(p,q),\dots,y_m(p,q))$. It is well known
that we have (this is called
the {\em addition formula})
\begin{equation*}
\tilde{\Pk}(p,q)=\frac{1}{d_{2\kappa}}\sum_{i=1}^{d_{2\kappa}}e_i(p)\overline{e_i(q)}.
\end{equation*}
As a consequence, from the expression
\begin{equation*}
[\Pk,\Pk]=\int_{\Gmn}\tilde{\Pk}(p,q)\tilde{\Pk}(q,p)dq
\end{equation*}
(i) follows. Moreover,
\begin{align*}
\sum_{p,q\in C}\tilde{\Pk}(p,q)&=\frac{1}{d_{2\kappa}}\sum_{p,q\in
  C}\big(\sum_{i=1}^{d_{2\kappa}}e_i(p)\overline{e_i(q)}\big)\\
&=\frac{1}{d_{2\kappa}}\sum_{i=1}^{d_{2\kappa}}\big(\sum_{p,q\in
  C}e_i(p)\overline{e_i(q)}\big)\\
&=\frac{1}{d_{2\kappa}}\sum_{i=1}^{d_{2\kappa}}\Big\vert\sum_{p\in
  C}e_i(p)\Big\vert^2\geq 0.
\end{align*}
hence (ii). More generally, for any function $\alpha:C\to \C$, we have:
\begin{align*}
\sum_{p,q\in C}\alpha(p)\overline{\alpha(q)}\tilde{\Pk}(p,q)=
&\frac{1}{d_{2\kappa}}\sum_{i=1}^{d_{2\kappa}}\Big\vert\sum_{p\in
  C}\alpha(p)e_i(p)\Big\vert^2\\
&\geq 0.
\end{align*}
Conversely, assume $F\in S_k$ is a polynomial with real coefficients,
such that, for any finite 
set $C\subset \Gmn$ and any function $\alpha:C\to \C$,
\begin{equation*}
\sum_{p,q\in C}\alpha(p)\overline{\alpha(q)}\tilde{F}(p,q)\geq 0,
\end{equation*}
and let us prove that $F$ expands on the $\Pk$ with non-negative coefficients.
Taking limits, we  have, for any $\alpha\in L^2(\Gmn)$,
\begin{equation*}
\int\int_{\Gmn}\alpha(p)\overline{\alpha(q)}\tilde{F}(p,q)dpdq
\geq 0
\end{equation*}
and hence, using the addition formula,
\begin{equation*}
\int\int_{\Gmn}\tilde{\Pk}(p,q)\tilde{F}(p,q)dpdq\geq 0.
\end{equation*}
If $F=\sum_{|\nu|\leq k}f_{\nu} \Pnu$, the left hand-side equals
$f_{\kappa}/d_{2\kappa}$, which proves that the coefficients
$f_{\kappa}$  are non-negative numbers. 
Using once again the addition formula, it is easy to show that
the product $\Pk\Pmu$ holds this general positivity property, and therefore
expands on the $\Pnu$ with non-negative coefficients.
%\eproof

\subsection{The principles of the LP bound}\label{subsec2.1}

The positivity property of the polynomials $\Pk$ is the basis of the
linear programming method to upper bound the cardinality of $\delta$-codes.

\medskip
\begin{definition} A Grassmannian code $C$ satisfying the 
constraint: 
\begin{equation*} \textrm{  For all } p\neq q\in C^2, d_c(p,q)\geq \delta.
\end{equation*}
is called a $\delta$-code.
\end{definition}
\medskip

\begin{proposition}\label{LPB} Assume  $F_{k}\in S_{k}$ satisfy:
\begin{enumerate}
\item $F_{k}=\sum_{|\kappa|\leq k} f_{\kappa}\Pk$ with $f_{\kappa}\geq 0$
for all $\kappa$, $f_0>0$
\item $F_{k}(y_1,\dots,y_m)\leq 0$ for all $(y_1,\dots,y_m)\in [0,1]^m$
such that \\$\sum_{i=1}^m y_i\leq m-\delta^2$
\end{enumerate}
Then, the following bound holds for the cardinality $|C|$ of any $\delta$-code:
\begin{equation*}
|C|\leq \frac{F_{k}(1,\dots,1)}{f_0}.
\end{equation*}
\end{proposition}
\medskip

\noindent\proof This is a standard argument, that we recall here.
Let $C$ be a $\delta$-code. As before, we let 
$\tilde{F_k}(p,q)=F_k(y_1(p,q),\dots,y_m(p,q))$. We calculate
\begin{equation*}
\sum_{p,q\in C}\tilde{F_k}(p,q)=\sum_{|\kappa|\leq k}f_{\kappa}\big(\sum_{p,q\in C}\tilde{\Pk}(p,q)\big)
\end{equation*}
Assumption (ii) leads to $\tilde{F_k}(p,q)\leq 0$ when $p\neq q$. The remaining terms
of the left hand-side, corresponding to $p=q$, contribute by
$|C|F_k(1,\dots,1)$. Assumption (i), together with the
positivity property of the polynomials $\Pk$ (Proposition \ref{PP} (ii)), show that all the terms 
of the right hand-side are non-negative. When $\kappa=(0)$, $\Pk=1$
and the contribution is $f_0|C|^2$. We obtain
\begin{equation*}
|C|F_k(1,\dots,1)\geq f_0|C|^2
\end{equation*}
equivalently
\begin{equation*}
|C|\leq \frac{F_k(1,\dots,1)}{f_0}.
\end{equation*}
It is worth noticing that equality in this inequality happens if and
only if, for all $1\leq |\kappa|\leq k$ such that $f_{\kappa}\neq 0$, $\sum_{p,q\in C}\tilde{\Pk}(p,q)=0$ 
and, for all $p\neq q\in C$, $\tilde{F_k}(p,q)=0$.
The first condition says that $C$ is a $2k$-design in the sense of
\cite{BCN} (when it holds for all $1\leq |\kappa|\leq k$),
and the second one that $C$ is an $F_k$-code in the sense of \cite{BBC}.
%\eproof

\subsection{The three-term relation and the Christoffel-Darboux formula}\label{subsec2.2}

We join here more material on the sequence of polynomials $\Pk$,
that will be of later use. The results presented here are essentially established in \cite{D}, except that we deal with symmetric polynomials.
Following \cite{D}, the (column) vector of the 
polynomials $\Pk$ with $|\kappa|=k$ is denoted by $\PP_k$. If necessary, we order the partitions of the same degree
in increasing lexicographic order.

We also set
\begin{equation*}
\sigma:=y_1+y_2+\dots+y_m
\end{equation*}
and, when necessary, we make the involved variables explicit, by writing
$\sigma(y)$ rather than $\sigma$.
The $\pi_k\times \pi_k$ diagonal matrix, denoted by  $D_k$, with entries
\begin{equation*}
D_k[\kappa,\kappa]=d_{2\kappa}:=\dim(V_n^{2\kappa})
\end{equation*}
is the  inverse of the Gram matrix of  $\PP_k$.

Next result is an analogue of the so-called ``three-term relation''.

\medskip
\begin{theorem}\label{3term}
For all $k\geq 1$, there exists matrices $A_k$, $B_k$, $C_k$,
of size respectively $\pi_k\times \pi_{k+1}$, $\pi_k\times \pi_{k}$, 
$\pi_k\times \pi_{k-1}$, such that:
\begin{equation*}
\sigma \PP_k=A_k\PP_{k+1}+B_k\PP_k+C_k\PP_{k-1}.
\end{equation*}
Moreover, $(D_kB_k)^t=D_kB_k$ and $D_kC_k=(D_{k-1}A_{k-1})^t$.
\end{theorem}
\medskip

\noindent\proof The polynomials $\sigma\Pk$ with $|\kappa|=k$
are symmetric of total degree $k+1$ so they afford a decomposition 
over the $(\Pmu)_{|\mu|\leq k+1}$. Moreover,
$[\sigma\Pk,\Pmu]=[\Pk,\sigma\Pmu]=0$ if $|\mu|\leq k-2$.

If $|\mu|=|\kappa|=k$, we have: $B_k[\kappa,\mu][\Pmu,\Pmu]=[\sigma\Pk,\Pmu]=
[\Pk,\sigma\Pmu]=B_k[\mu,\kappa][\Pk,\Pk]$, which proves that the matrix 
$D_kB_k$ is symmetric. The same argument shows that 
$D_kC_k=(D_{k-1}A_{k-1})^t$.

%\eproof

\noindent{\bf Notations:} We want to define  $\ki$
(respectively $\kii$) to be the partition 
obtained from $\kappa$ by increasing (respectively decreasing) 
the $i$-th part $\kappa_i$ by one. This is not possible for all $i$,
since the result should be also a partition, i.e. the new parts
should be in decreasing order. Hence we define (where $\kappa_{m+1}:=0$):
\begin{equation*}
\begin{cases}
u(\kappa):=\{1\}\cup \{i\in [2,m]\mid \kappa_{i-1}>\kappa_i\}\\
d(\kappa):=\{i\in [1,m]\mid \kappa_{i}>\kappa_{i+1}\}
\end{cases}
\end{equation*}
The set $u(\kappa)$ is the set of indices $i$ for which $\ki$ makes sen\-se
(respectively $d(\kappa)$ for $\kii$).
Moreover, if $|\kappa|=k$,
$|\ki|=k+1$ and $|\kii|=k-1$.

Otherwise explicitly mentioned, in the rest of this paper,
$\kappa$, $\kappa'$ are partitions
of degree $k$, while $\mu$, $\mu'$ are partitions of degree $k+1$ and
$\nu$, $\nu'$ are partitions of degree $k-1$.

\smallskip

\medskip
\begin{proposition}\label{coef} The following properties hold:
\begin{enumerate}
\item For all $\kappa$, $\mu$ and $\kappa'$, $A_k[\kappa,\mu]\geq 0$, and $
B_k[\kappa,\kappa']\geq 0$.
\item The coefficients of the matrix $A_k$ are equal to
zero, except the coefficients $A_k[\kappa,\ki]$, which are positive. 
\end{enumerate}
\end{proposition}
\medskip

\noindent\proof The first assertions are equivalent to:
$[\sigma \Pk,\Pmu]\geq 0$ for all $\ka$, $\mu$ of any degree. But  
$[\sigma \Pk,\Pmu]=[1,\sigma\Pk\Pmu]$ and 
$\sigma=m(1-\frac{m}{n})P_1+\frac{m^2}{n}$. Joint with 
Proposition \ref{PP}(iii), we obtain $[1,\sigma\Pk\Pmu]\geq 0$.

The coefficients $A_k[\ka,\mu]$ can be more precisely calculated,
using (\cite[Lemma 7.5.7]{M}). Since we do not normalize the
polynomials $\Ck$ in the same way, we introduce coefficients
$\bb{\mu}{\ka}$ such that
\begin{equation*}
\sigma C_{\kappa}=\sum_{|\mu|=k+1}\bb{\mu}{\ka}\Cmu.
\end{equation*}
They differ by a positive multiplicative factor from the generalized binomial
coefficients $\binom{\mu}{\kappa}$
defined in \cite{M};
see also \cite{La1}. Then we have
\begin{equation*}
A_k[\kappa,\mu]=\bb{\mu}{\ka}\left(\frac{\beta_{\mu}}{\beta_{\ka}}\right)^{-1}.
\end{equation*}
It is known that the generalized binomial
coefficients $\binom{\mu}{\kappa}$ are equal to zero when
$\mu$ is not equal to one of the $\ki$; consequently the same holds for $A_k[\kappa,\mu]$.
Moreover, since $\binom{\ki}{\ka}>0$, also $\bb{\ki}{\ka}>0$ and $A_k[\kappa,\ki]\neq 0$.
%\eproof

\medskip
\begin{theorem}[Christoffel-Darboux Formula]\label{CDF}
Let
\begin{equation*}
\Qk:=\sum_{|\mu|=k+1}A_k[\kappa,\mu]\Pmu\in \Pi_{k+1}.
\end{equation*}
With the previous notations, we have:
\begin{enumerate}
\item For all $k\geq 0$,
\begin{align*}
&\sum_{|\nu|\leq k} d_{2\nu}\Pnu(x)\Pnu(y)=\\
&\frac{\sum_{\substack{|\kappa|=k}}d_{2\kappa}\big(\Qk(x)\Pk(y)-\Pk(x)\Qk(y)\big)}
{\sigma(x)-\sigma(y)}
\end{align*}

\item Moreover, if $\epsilon:=\sum_{i=1}^m \frac{\partial}{\partial y_i}$,

\begin{align*}
&\sum_{|\nu|\leq k} d_{2\nu}\big(\Pnu(y)\big)^2=\\
&\sum_{\substack{|\kappa|=k}}\frac{d_{2\kappa}}{m}\Big(\big(\epsilon\Qk(y)\big)\Pk(y)-
\big(\epsilon\Pk(y)\big)\Qk(y)\Big).
\end{align*}
\end{enumerate}
\end{theorem}
\medskip

\noindent\proof 
The proof of (i) is the same as \cite{D}[Theorem 3.5.3]. 
Note that we cannot hope for a formula
for each $x_i$ like in \cite{D}, since we should stick to 
symmetric polynomials. If

\begin{align*}
\Sigma_s:&=(A_s\PP_{s+1}(x))^tD_s\PP_s(y)-\PP_s(x)^tD_sA_s\PP_{s+1}(y)\\
&=\sum_{\substack{|\kappa|=s}}d_{2\kappa}\big(\Qk(x)\Pk(y)-\Pk(x)\Qk(y)\big),
\end{align*}
from the ``three-term relation'' of Theorem \ref{3term}, we have:

\begin{equation*}
\Sigma_s-\Sigma_{s-1}=(\sigma(x)-\sigma(y))\PP_s(x)^tD_s\PP_s(y).
\end{equation*}
The formula (i) follows from summing up these identities, for $1\leq
s\leq k$.

\smallskip
In the equation (i), we replace $\Qk(x)\Pk(y)-\Pk(x)\Qk(y)$ by

\begin{equation*}
\Qk(x)\big(\Pk(y)-\Pk(x)\big)-\Pk(x)\big(\Qk(y)-\Qk(x)\big).
\end{equation*}
Then, if we specialize
$x_2=y_2, \dots,x_m=y_m$ and let $x_1$ tend to $y_1$, we obtain 

\begin{align*}
&\sum_{|\nu|\leq k} d_{2\nu}\big(\Pnu(y)\big)^2=\\
&\sum_{\substack{|\kappa|=k\\|\mu|=k+1}}d_{2\kappa}A_k[\kappa,\mu]
\Big(\big(\frac{\partial\Qk}{\partial y_1}(y)\big)\Pk(y)-\big(\frac{\partial\Pk}{\partial y_1}(y)\big)\Qk(y)\Big).
\end{align*}

The same identity holds when one replaces $y_1$ by any $y_i$;
if we sum up all these identities, we obtain the more symmetric
formula (ii).

%\eproof

\begin{remark} 
The left hand side of the Christoffel-Darboux formula

\begin{equation*}
K_k(x,y):=\sum_{|\nu|\leq k} d_{2\nu}\Pnu(x)\Pnu(y).
\end{equation*}
is the {\em reproducing kernel} of the space of symmetric polynomials
of degree at most $k$. It satisfies
the characteristic property: for all $Q\in S_k$,  $[K_k(x,.),Q]=Q(x)$.

\end{remark}

\subsection{An LP bound from Christoffel-Darboux formula}\label{subsec2.3}

In the classical cases, Christoffel-Darboux formula is involved 
in the setting up of bounds of the type
$|C|\leq M(\delta)$ where $M(\delta)$ is an explicit function of $\delta$.
Usually the running interval of $\delta$ is divided into subintervals, related to the zeros 
of the zonal polynomials. This is the line followed 
in \cite{MRRW}, and also in \cite {KL}; see \cite{L} for a unified
presentation.
In this section, we follow this method, and analyze the difficulties 
arising from  the several variables situation.

The numerator, of degree $k+1$,
 of the right hand side 
of Christoffel-Darboux formula (Theorem \ref{CDF}(i)) is denoted by  $N_{k+1}(x,y)$. We consider 
the polynomial in the variables $y_1,\dots ,y_m$, of degree $2k+1$,

\begin{equation*}
\begin{split}
F_{2k+1}(x,y):&=-N_{k+1}(x,y)K_k(x,y)\\
&=\frac{N_{k+1}(x,y)^2}{\sigma(y)-\sigma(x)}.
\end{split}
\end{equation*} 
In order to make clear that only the $y_1,\dots ,y_m$ are variables,
while the $x_1,\dots ,x_m$ will specialize to real values, we denote
it by $F_{2k+1}(x,\cdot)$.

\medskip
\begin{proposition}\label{pp1} Let $x\in [0,1]^m$ satisfy $\sigma(x)>s:=m-\delta^2$.
Assume the following conditions hold:
\begin{enumerate}
\item For all $\kappa$, $|\kappa|\leq k$, $\Pk(x)\geq 0$
\item For all $\kappa$, $|\kappa|= k$, $\Qk(x)\leq 0$
\end{enumerate}
Then, $F_{2k+1}(x,\cdot)$ satisfies the conditions required in Proposition \ref{LPB}.
\end{proposition}
\medskip

\noindent\proof We have:
\begin{equation*} 
F_{2k+1}(x,y)=\frac{N_{k+1}(x,y)^2}{\sigma(y)-\sigma(x)}
\end{equation*}
hence condition (ii) is satisfied when $s<\sigma(x)$.

To prove condition (i), we point out that,
if $F$ and $G$ are two polynomials with non-negative
 coefficients
on the $\Pk$, then the product $FG$ holds the same property.
This is a direct consequence of Proposition \ref{PP}(iii).

From the definition of $K_k(x,y)$, its coefficient on $\Pnu$
with $|\nu|\leq k$ equals $d_{2\nu}\Pnu(x)$ (and for higher degree
partitions it is zero).
On the other hand, 
\begin{align*}
-N_{k+1}(x,y)&=\sum_{|\kappa|=k} d_{2\kappa}\big(\Pk(x)\Qk(y)- \Qk(x)\Pk(y)\big)\\
&=\sum_{|\mu|=k+1} \Big(\sum_{|\kappa|=k}d_{2\kappa}A_k[\kappa,\mu]\Pk(x)\Big)\Pmu(y)\\
&\quad -\sum_{|\kappa|=k}d_{2\kappa}\Qk(x)\Pk(y)
\end{align*}

The coefficient $A_k[\kappa,\mu]$ is always non-negative. 
Clearly, under the conditions of the proposition,
the coefficients of $-N_{k+1}(x,y)$ on the $\Pk$ and $\Pmu$ 
are non-negative.

%\eproof

\medskip
\begin{corollary}  Assume $x$ satisfies the conditions of Proposition \ref{pp1}.
Then, for all $\delta$-code $C$, 
\begin{equation*}
|C| \leq 
\frac{\big(m-\s(x)\big)\big(\sum_{|\nu|\leq k}d_{2\nu}\Pnu(x)\big)^2}
{-\sum_{|\kappa|=k}d_{2\kappa}\Pk(x)\Qk(x)}
\end{equation*}
\end{corollary}
\medskip

\noindent\proof In order to apply Proposition \ref{LPB},
we are left with the computation of $f_0$ and of $F_{2k+1}(x,(1,\dots,1))$.
Since $F_{2k+1}(x,y)=K_k(x,y)^2(\sigma(y)-\sigma(x))$, we have
\begin{equation*}
F_{2k+1}(x,(1,\dots,1))=\big(\sum_{|\nu|\leq k}d_{2\nu}\Pnu(x)\big)^2(m-\sigma(x)).
\end{equation*}
Using the orthogonality of the $\Pk$, we obtain
\begin{align*}
f_0:&=[F_{2k+1}(x,\cdot),1]=-[K_k(x,\cdot),N_{k+1}(x,\cdot)]\\
&=-[\sum_{|\nu|\leq k}d_{2\nu}\Pnu(x)\Pnu,\sum_{|\kappa|=k}d_{2\kappa}\Qk(x)\Pk]\\
&=-\sum_{|\kappa|=k}d_{2\kappa}\Pk(x)\Qk(x).
\end{align*}
%\eproof

The main problem with this approach, is that, in general, we don't even know
if the inequalities (i) and (ii) of Proposition \ref{pp1} have a solution $x$.
In case these inequalities define a non empty area of $\R^m$, a second
problem would be to optimize the choice of $x$ in this area.
In the classical case $m=1$, $Q_k=P_{k+1}$ (up to a positive multiplicative factor). The interlacing property of the real zeros of the orthogonal polynomials
$P_k$, ensures that one can take $x\in [z_k,z_{k+1}]$, where $z_k$ is the
largest zero of $P_k$, so that $P_{k+1}(x)\leq 0$ and $P_i(x)\geq 0$ for all
$i\leq k$. Moreover, one uses asymptotic estimates of these zeros
to derive an asymptotic bound for the size of codes. 

In the general case $m\geq 2$, we don't have such tools to deal with 
the inequalities of Proposition \ref{pp1}, which seem to be intractable
in general. The first case $k=1$, leading to a polynomial of degree $3$,
is however discussed in the next section.
On the other hand, one can think of the zeros of orthogonal polynomials
in one variable as being the eigenvalues of the so-called
Jacobi matrices associated to the sequence of polynomials. We 
study in section \ref{sec4} the eigenvalues of the analogous 
matrices in the general case,
and derive bounds for codes, which contain as a special case the
bound obtained from a possible solution of these inequalities.

\section{LP bounds of small degree}\label{sec3}

We take the following notations: let $s:=m-\delta^2$, the maximal
value of $\sigma$ among pairs of points of a code $C$. We are looking for 
a function $M(s)$ such that $|C|\leq M(s)$. Obviously, $M(s)$ is an increasing
function. In this section, we discuss the cases of small degree $k$,
trying to optimize the choice of $F_k$ in Proposition \ref{LPB}

\subsection{Degree $1$}\label{subsec3.1}

Let $F_1=1+f_1P_1$, with $f_1\geq 0$ (condition (i)). We have $P_1=\frac{n}{m(n-m)}(\sigma-\frac{m^2}{n})$. 

When $\s\in [0,s]$, $1+f_1P_1$ should be non-positive (condition (ii)). 
Therefore, The zero
of $1+f_1P_1$ should be greater than $s$. It leads to the condition:

\begin{equation*}
s-\frac{m^2}{n}\leq -\frac{m(n-m)}{nf_1 }.
\end{equation*}

Since $f_1\geq 0$, we obtain the necessary condition
$s\leq \frac{m^2}{n}$. The smallest value for $f_1$ is then

\begin{equation*}
f_1=\frac{-m(1-m/n)}{s-\frac{m^2}{n}}
\end{equation*}
corresponding to a polynomial proportional to $\s-s$. 
We obtain the bound

\begin{equation*}
\textrm{if }s< \frac{m^2}{n},\  |C|\leq \frac{m-s}{\frac{m^2}{n}-s}
\end{equation*}
which is the so-called simplex bound proved in 
\cite{CHS}.

\subsection{Degree $2$}\label{subsec3.2}

We restrict ourselves to polynomials which are divisible by $\s-s$.
Then, such polynomials are polynomials in $\s$. We write:

\begin{equation*}
F_2=(\s-s)(\s-b)=f_2P_2+f_{11}P_{11}+f_1P_1+f_0.
\end{equation*}
with the condition that $b\leq 0$. With $t=s-m^2/n$, we find:

\begin{equation*}
\begin{split}
f_2:=& \frac{m(m+2)(n-m)(n-m+2)}{3(n+2)(n+4)}\\
f_{11}:=&\frac{2m(m-1)(n-m)(n-m-1)}{3(n-2)(n-1)}\\
f_1:=&m\left(1-\frac{m}{n}\right)\left(\frac{m^2}{n}+\frac{4(n-2m)^2}{n(n-2)(n+4)}-t -b\right)\\
f_0:=&\frac{2m^2(n-m)^2}{n^2(n-1)(n+2)}-\frac{m^2}{n}t+bt
\end{split}
\end{equation*}

The condition $f_0>0$, when $t>0$,  is equivalent  to 

\begin{equation}\label{e3}
b> \frac{m^2}{n}\left(1-\frac{2(n-m)^2}{tn(n-1)(n+2)}\right)
\end{equation}
(and when $t<=0$ is always fulfilled), which implies

\begin{equation}\label{e4}
t< \frac{2(n-m)^2}{n(n-1)(n+2)}.
\end{equation}

The condition $f_1\geq 0$ is equivalent to 

\begin{equation}\label{e5}
b \leq \frac{m^2}{n}+\frac{4(n-2m)^2}{n(n-2)(n+4)}-t.
\end{equation}

One can check that the right hand side of (\ref{e5}) is positive 
for $m\geq 2$, when $t$ satisfies (\ref{e4}).

The bound $B=(f_2+f_{11}+f_1+f_0)/f_0$ equals

\begin{equation}\label{e6}
B=\delta^2\frac{m-b}{f_0}.
\end{equation}

As a function of $b$, it is decreasing when 
$t\in [\frac{-2m(n-m)}{n(n-1)(n+2)}, \frac{2(n-m)^2}{n(n-1)(n+2)}[$, and hence
the best choice of $b$ is $b=0$. We obtain the bound:

\medskip
\begin{theorem} If $s\in \left]0, \frac{m^2}{n}+\frac{2(n-m)^2}{n(n-1)(n+2)}\right[$,
\begin{equation}\label{B2}
|C|\leq \frac{n}{m}\Big(\frac{m-s}{-s+\frac{m^2}{n}+\frac{2(n-m)^2}{n(n-1)(n+2)}}\Big)
\end{equation}
\end{theorem}
\medskip

This bound, which is an increasing function of s,
improves  on the simplex bound when $s\geq
\frac{m^2}{n}-\frac{2m(n-m)}{n(n-1)(n+2)}$. Their common value at $s=
\frac{m^2}{n}-\frac{2m(n-m)}{n(n-1)(n+2)}$ is $\binom{n+1}{2}$.
However, the orthoplex bound proved in \cite[(5.6)]{CHS} 
reads: 
\begin{align*}
s<m^2/n &\Rightarrow |C|\leq \binom{n+1}{2}\\
s=m^2/n &\Rightarrow |C|\leq (n-1)(n+2)
\end{align*}
and is better
than (\ref{B2}) in the range 
$]\frac{m^2}{n}-\frac{2m(n-m)}{n(n-1)(n+2)}, \frac{m^2}{n}[$.
If we plug in (\ref{B2}) the value $s=m^2/n$, we find that
 $|C|\leq \frac{(n-1)(n+2)}{2(1-m/n)}$  which is better than the
orthoplex bound when $m<n/2$. 
We recall that the orthoplex bound is attained 
for a family of codes
with $n=2^i$, $m=n/2$, constructed in \cite[Theorem 1]{CHRSS}. These
codes are also optimal $6$-designs (see \cite{B}).

\subsection{Degree $3$}\label{subsec3.3}

We do not study general polynomials of degree $3$ but rather 
apply the approach described in subsection \ref{subsec2.3}

The polynomial $F_3$ has degree $3$, and is again a polynomial in 
$\s$. In the following, we calculate the best choice for
$x$ (and discuss its existence).  Let 
$u:=\s(x)-m^2/n$. We should have:

\begin{enumerate}

\item $u\geq s-m^2/n$
\item $u\geq 0$  (Condition  (i))
\item $u^2-\frac{4(n-2m)^2}{n(n-2)(n+4)}u-\frac{2m^2(n-m)^2}{n^2(n-1)(n+2)}\leq 0$ (Condition  (ii))
\end{enumerate}

The polynomial of degree $2$ occurring in (iii) has a positive discriminant,
and a unique positive root that we shall denote by $u_2$. 
Let $b$ and $c$ be the coefficients of this polynomial, so that it is
equal to $u^2-bu-c$, and let $d:=\frac{2m(n-m)}{n(n-1)(n+2)}$.
The bound is then equal to:

\begin{equation*}
B(u):=-\frac{(n-1)(n+2)(u+d)^2(m-u-m^2/n)}{2u(u^2-bu-c)}.
\end{equation*}

The calculation of $B'(u)$ shows that it is increasing in the range
$[u_1,u_2]$ (the numerator has the form:
$u+d$ times a degree $3$ polynomial with a unique real root $u_1$). 
Hence, for $s\in [u_1+\frac{m^2}{n},u_2+\frac{m^2}{n}]$,
the best choice for $u$ is $u=s-m^2/n$.

We obtain:

\medskip
\begin{theorem} Let $d=\frac{2m(n-m)}{n(n-1)(n+2)}$,
$b=\frac{4(n-2m)^2}{n(n-2)(n+4)}$,
$c=\frac{2m^2(n-m)^2}{n^2(n-1)(n+2)}$.

\begin{equation*}
\begin{split}
\textrm{If }&s\in \left]\frac{m^2}{n}, \frac{m^2}{n}+\frac{b}{2}+
\sqrt{\frac{b^2}{4}+c}\right[,\\
&|C|\leq \frac{(m-s)(s-\frac{m^2}{n}+d)^2(n-1)(n+2)}
              {2(s-\frac{m^2}{n})(-(s-\frac{m^2}{n})^2
               +b(s-\frac{m^2}{n})+c)}
\end{split}
\end{equation*}
\end{theorem}
\medskip

\section{The endomorphisms $T_k$}\label{sec4}

We introduce an endomorphism $T_k: S_k\to S_k$ which eigenvalues will
play the role of the zeros of the zonal polynomials in the rank one case.

\medskip
\begin{proposition}\label{p61} Let 
\begin{align*}
T_k: S_k & \to S_k\\
P&\mapsto \pr_{S_k}(\sigma P)
\end{align*}
where the orthogonal 
projection on $S_k$ is denoted by $\pr_{S_k}$ (note that, in general, $\sigma P$ does not belong to $S_k$
but rather to $S_{k+1}$). 

The endomorphism $T_k$ is a symmetric endomorphism of $S_k$,
and is an isomorphism.
\end{proposition}
\medskip

\noindent\proof We have, for all $P,Q\in S_k$, 
$[T_k(P),Q]=[\sigma P,Q]=[P,\sigma Q]=[P, T_k(Q)]$.
Moreover, $[\sigma P,P]=[\sigma, P^2]>0$ unless $P=0$, because of the
positivity of the measure on $\R[y_1,\dots,y_m]^{S_m}$. Thus $T_k$ is injective.

%\eproof

Let $J_k$ be the matrix of this endomorphism in the basis 
$\{\Pk, |\kappa|\leq k\}$.
From the three-term relation (Theorem \ref{3term}), $J_k$ is the block-tridiagonal matrix: 

\begin{equation}\label{jk}
J_k=\begin{pmatrix}B_0 & A_0& &&   &\\
                   C_1&B_1&A_1 &&&  \\
                    &C_2&B_2&A_2&&\\
                     &&C_3&\ddots&\ddots&\\
                       &&&\ddots&\ddots&A_{k-1}\\

                       & &&&C_k&B_{k}
\end{pmatrix}
\end{equation}

It is worth noticing that the matrix $J_k$ itself is not symmetric, 
because the polynomials $\Pk$ are not of norm $1$. We shall later
introduce and calculate the symmetric matrix ${J'}_k$ obtained in the
normalized basis. 
\smallskip

In the end, we shall need some very precise information on the
coefficients of ${J'}_k$. For the moment, the only, but crucial, property
that we will exploit is the fact that it is {\em non-negative and
  irreducible}.

\smallskip

\medskip
\begin{lemma}\label{eigenv} The eigenvalues of $T_k$ are real, and belong to $]0,m[$. 
The maximal eigenvalue of $T_k$, denoted by $\lambda_k$,
is of multiplicity $1$, and possesses an eigenvector with positive 
coordinates. Moreover, $\lambda_{k-1}<\lambda_k$.

\end{lemma}
\medskip

\noindent\proof The matrix $J_k$ is non-negative and irreducible in the sense of 
\cite{G}, because of Proposition \ref{coef} (note that
the coefficients $A_k[\kappa,\ki]$ are positive).
Moreover, it is the matrix of a symmetric endomorphism, so its eigenvalues are real. From \cite[Perron-Frobenius Theorem]{G}, it follows that 
the maximal eigenvalue has 
multiplicity equal to $1$, and that, if $v$ is an eigenvector, either $v$
or $-v$ has positive coordinates.
Let us now prove that all its eigenvalues belong to $]0,m[$.

For any $v\in S_k$, $v\neq 0$,
we have $[\sigma v,v]=\int \sigma v^2 d\mu(y)$, where $d\mu$ is a
positive measure. 
We integrate on the domain $[0,1]^m$, on which $0\leq \sigma\leq m$,
hence $0<[\sigma v,v]<m[v,v]$.
If $v$ is an eigenvector of $T_k$ associated with an eigenvalue $\lambda$,
we have $[\sigma v,v]=[\lambda v,v]=\lambda[v,v]$, so we can conclude that $0<\lambda<m$.

Now let $v$ be an eigenvector of $T_{k-1}$  for $\lambda_{k-1}$, assumed to be of norm $1$. We have 

$$\sigma v=\lambda_{k-1} v +u$$
with $u\in \Pi_k$. Obviously, since $\deg(\sigma v)=1+\deg(v)$, $v$
must be of degree exactly $k-1$ (and $u\neq 0$).

Since 

$$\lambda_k=\max_{x\in S_k\setminus\{0\}}\frac{[T_k(x),x]}{[x,x]},$$
we have
$[T_k(v),v]\leq \lambda_k$. But $[T_k(v),v]=[\sigma v,v]=\lambda_{k-1}$.
The equality $\lambda_{k-1}=\lambda_k$ would mean that $v$ is an eigenvector
of $T_k$, which is not possible since it has degree $k-1$.

%\eproof

In the case $m=1$, the eigenvalues of $T_k$ are exactly 
the zeros of the polynomial
$P_{k+1}$. In the general case, we prove in next lemma that common
zeros of the polynomials $\Qk$ give some of the eigenvalues. However,
we do not know if such common zeros do exist, neither 
if all of the eigenvalues  are obtained that way (and may be
it is not so important):

\medskip
\begin{lemma}
Let $\alpha\in [0,1]^m$ be a common zero of the polynomials
\begin{equation*}
Q_{\kappa}:=\sum_{|\mu|=k+1} A_k[\kappa,\mu]\Pmu,
\end{equation*} 
for all $\kappa$, 
$|\kappa|=k$. Then, $v:=\sum_{|\nu|\leq k}d_{2\nu}\Pnu(\alpha)\Pnu$ 
is an eigenvector of $T_k$ for the eigenvalue $\sigma(\alpha)$.

\end{lemma}
\medskip

\noindent\proof It is immediate from Christoffel-Darboux formula (Theorem
\ref{CDF}(i)). 
If $Q_{\kappa}(\alpha)=0$ for all $\kappa$, $|\kappa|=k$, we have

\begin{equation*}
(\sigma(\alpha)-\sigma(y))v=-\sum_{\substack{|\kappa|=k\\|\mu|=k+1}}d_{2\kappa}A_k[\ka,\mu]\Pk(\alpha)\Pmu(y)\in \Pi_{k+1}
\end{equation*}
and, therefrom,

\begin{equation*}
\sigma(\alpha)v=T_k(v).
\end{equation*}

%\eproof

We now show how to obtain a bound for the size of $\delta$-codes, as a
function of $\delta$. Therefore, in order to cope with any possible
$\delta$, we must perturb the endomorphism $T_k$ as explained next:

\medskip
\begin{theorem}\label{Teps} Let $\epsilon\in \R^{\pi_k}$, 
with $\epsilon_{\kappa}\geq 0$.
Let $T_k^{\epsilon}$
be the endomorphism defined on $S_k$ by

\begin{equation*}
T_k^{\epsilon}(v)=T_k(v)-\epsilon* v_k
\end{equation*}
where $\epsilon* v_k:=
\sum_{|\kappa|=k}\epsilon_{\kappa}v_{\kappa}P_{\kappa}$. 

\begin{enumerate}
\item $T_k^{\epsilon}$ has a unique maximal eigenvalue 
$\lambda_k^{\epsilon}$, of multiplicity one, possessing an eigenvector
$v^{\epsilon}$ with positive coefficients. Moreover, if $\epsilon\neq 0$,

\begin{equation*}
\lambda_{k-1}<\lambda_k^{\epsilon}< \lambda_k
\end{equation*}

\item Let $\epsilon \neq 0$. Any $\delta$-code $C$ such that $s=m-\delta^2< \lep$ 
satisfies

\begin{equation*}
|C|\leq
\frac{\big(\sum_{|\kappa|=k}v_{\kappa}^{\epsilon}(\epsilon_{\kappa}+a_{\kappa})\big)^2}
{(m-\lep)\big(\sum_{|\kappa|=k}d_{2\kappa}^{-1}\epsilon_{\kappa}{v_{\kappa}^{\epsilon}}^2\big)}
\end{equation*}
where $a_{\kappa}:=Q_{\kappa}(1,\dots,1)=\sum_{|\mu|=k+1}A_k[\kappa,\mu]$.
\end{enumerate}
\end{theorem}
\medskip

\noindent\proof

(i) The matrix $J_k^{\epsilon}$ of $T_k^{\epsilon}$ is equal to $J_k$, 
except the diagonal
elements lying in $B_k$. Replacing $J_k^{\epsilon}$ by $J_k^{\epsilon}+M\Id$
for some appropriate $M$, we obtain a non-negative matrix which is irreducible
so its largest eigenvalue has multiplicity one and has an associated 
eigenvector  
with positive coordinates. It remains true
for $J_k^{\epsilon}$. Since, when $\epsilon \neq 0$,
 $J_k^{\epsilon}<J_k$, we have $\lep<\la_k$.
The proof of the inequality $\la_{k-1}<\lep$ is the same as the one of
$\la_{k-1}<\la_k$.

(ii) We have $\sigma \ve=\lep \ve+\epsilon * v_k^{\epsilon} +u$ 
where $u\in \Pi_{k+1}$. We need to compute $u$, and we set 
$u=\sum_{|\mu|= k+1}u_{\mu}\Pmu$.
Let $\mu$, $|\mu|=k+1$, we have:
\begin{align*}
u_{\mu}[\Pmu,\Pmu]&=[u,\Pmu]=[\sigma \ve,\Pmu]\\
&=\sum_{|\kappa|\leq k} v_{\kappa}^\epsilon[\sigma \Pk, \Pmu]\\
&=\sum_{|\kappa|= k} v_{\kappa}^\epsilon A_k[\kappa,\mu][\Pmu,\Pmu]\\
\end{align*}
and we obtain $u_{\mu}=\sum_{|\kappa|= k} v_{\kappa}^\epsilon
A_k[\kappa,\mu]$.
We have found
\begin{equation*}
u=\sum_{|\mu|=k+1}\big(\sum_{|\kappa|=k}v_{\kappa}^{\epsilon}A_k[\kappa,\mu]\big)\Pmu
=\sum_{|\kappa|=k}v_{\kappa}^{\epsilon}Q_{\kappa},
\end{equation*}
hence the {\em ``generalized Christoffel-Darboux formula''}:
\begin{equation}\label{cdf3}
\ve=\sum_{s=0}^k v_s^{\e} \PP_s=
\frac{\sum_{|\kappa|=k}v_{\kappa}^{\epsilon}(\epsilon_{\kappa}P_{\kappa}+Q_{\kappa})}
{\sigma-\lep}.
\end{equation}

Now we proceed like in Proposition \ref{pp1}.
Let the numerator of the right hand side be denoted by  $N_{k+1}(y)$, and
let 
\begin{equation*}
F_{2k+1}(y):=\frac{N_{k+1}(y)^2}{\sigma(y)-\lep}=N_{k+1}(y)\ve.
\end{equation*}
We have:
\begin{align*}
f_0=[F_{2k+1},1]&=[N_{k+1},\ve]\\
&=[\sum_{|\kappa|=k}v_{\kappa}^{\epsilon}\epsilon_{\kappa}P_{\kappa},\sum_{|\kappa|=k}v_{\kappa}^{\epsilon}P_{\kappa}]\\
&=\sum_{|\kappa|=k}{v_{\kappa}^{\epsilon}}^2\epsilon_{\kappa}d_{2\kappa}^{-1}.
\end{align*}
Since the coefficients of $\epsilon$ and of $\ve$ are non-negative  
numbers, and $f_0\neq 0$ when $\epsilon \neq 0$, it follows that $F_{2k+1}$
satisfies the condition (i) of Proposition \ref{LPB}. Condition (ii)
is clearly fulfilled if $s<\lep$. We calculate
\begin{equation*}
F_{2k+1}(1,\dots,1)=\frac{\big(\sum_{|\kappa|=k}v_{\kappa}^{\epsilon}(\epsilon_{\kappa}+a_{\kappa})\big)^2}{m-\lep}.
\end{equation*}
hence the announced bound.

%\eproof

Let us show that  we have indeed generalized the situation
described in subsection \ref{subsec2.3} and Proposition \ref{pp1}.
Let $x\in \R^m$ such that $\Pk(x)>0$,
$Q_{\kappa}(x)\leq 0$ for all $|\kappa|=k$, and 
$\Pk(x)\geq 0$ for all $|\kappa|\leq k$.
Let $\epsilon\in \R^{\pi_k}$  be defined by: 
$\epsilon_{\kappa}=-\Qk(x)/\Pk(x)$. We can show that $\lep=\sigma(x)$.
Indeed, from (i) and (ii) of the proposition, 

\begin{equation*}
\ve(y)=
\frac{\sum_{|\kappa|=k}v_{\kappa}^{\epsilon}\big(-(\Qk(x)/\Pk(x))P_{\kappa}(y)+Q_{\kappa}(y)\big)}
{\sigma(y)-\lep}.
\end{equation*}

When we let $y$ tend to $x$, the numerator tends to $0$. Since the coordinates
of $\ve$ are positive and  $\Pk(x)\geq 0$ for all $|\kappa|\leq k$,
the left hand side cannot be equal to zero when $y=x$ ($P_0=1$). 
So the denominator 
also tends to zero, and $\lep=\sigma(x)$. The Christoffel-Darboux formula
(Theorem \ref{CDF}(i)) shows that $\ve=\sum_{|\kappa|\leq k} d_{2\kappa}\Pk(x)\Pk$.

When $m=1$, $\pi_k=1$ and any $\epsilon\geq 0$ is of this form. When $m\geq 2$,
it is not clear.. It is not even clear that at least one $x$ satisfying
these inequalities exists.

\medskip

Another natural question concerns the values that
$\lambda_k^{\epsilon}$ takes. It is hoped of course that all values in
the interval $]\lambda_{k-1},\lambda_k]$ are attained.
We have defined a mapping from 
$[0,+\infty[^{\pi_k}$ to $]\lambda_{k-1},\lambda_k]$, sending $\epsilon$
to $\lep$, which is continuous, hence the image in an interval,
containing $\lambda_k$, since clearly it is the image of $\epsilon=0$.
Let us prove that $\lep$ tends to $\lambda_{k-1}$ when 
$\epsilon$ tends to $+\infty$. To that end, we use the following
inequality,
valid for any non-negative matrix J with maximal eigenvalue $\lambda$
(\cite{G}):
\begin{align*}
\textrm{For all }x, x_i>0,\ \lambda\leq\sup_{i}\frac{(xJ)_i}{x_i}.
\end{align*}
This inequality remains true for the matrix $J_k^{\epsilon}$, although it is  not
non-negative, because we can apply it to some  $J_k^{\epsilon}+M\Id$,
an argument that we have already called for.
We choose for $x\in \R^{s_k}$ a vector,
which first $s_{k-1}$ coefficients constitute a positive eigenvector
of $J_{k-1}$ for the eigenvalue $\lambda_{k-1}$. 
Its last $\pi_k$ coordinates are denoted by
$u=(u_{\kappa})_{|\kappa|=k}$. We have:

\begin{equation*}
\begin{cases}
\textrm{If } |\nu|\leq k-2,
\frac{(xJ_k^{\epsilon})_{\nu}}{x_{\nu}}=\lambda_{k-1}\\
\textrm{If } |\nu|= k-1,
\frac{(xJ_k^{\epsilon})_{\nu}}{x_{\nu}}=\lambda_{k-1}
+\frac{\sum_{|\kappa|=k}u_{\kappa}C_{k}[\kappa,\nu]}{x_{\nu}}\\
\textrm{If } |\kappa|= k,
\frac{(xJ_k^{\epsilon})_{\kappa}}{x_{\kappa}}=
\frac{\sum_{|\nu|=k-1}x_{\nu}A_{k-1}[\nu,\kappa]}{u_{\kappa}}\\
\qquad\qquad\qquad\qquad\qquad+(B_k[\kappa,\kappa]-\epsilon_{\kappa}).
\end{cases}
\end{equation*}

The last equality relies on a result that is only proved in Section
5, Proposition \ref{ppp1}(i), namely that $B_k[\kappa,\kappa']=0$ when
$\kappa\neq \kappa'$.

Let us now choose an arbitrary small $\alpha>0$; we can choose the
coefficients $u_{\kappa}>0$ such that 
$\frac{\sum_{|\kappa|=k}u_{\kappa}C_{k}[\kappa,\nu]}{x_{\nu}}\leq
\alpha$ for all $\nu$ with $|\nu|=k-1$. Then we can choose
$\epsilon_{\kappa}>0$ such that 
$\frac{\sum_{|\nu|=k-1}x_{\nu}A_{k-1}[\nu,\kappa]}{u_{\kappa}}+(B_k[\kappa,\kappa]-\epsilon_{\kappa})=0$.
We are left with:
\begin{equation*}
\begin{cases}
\textrm{If } |\nu|\leq k-2,
\frac{(xJ_k^{\epsilon})_{\nu}}{x_{\nu}}=\lambda_{k-1}\\
\textrm{If } |\nu|= k-1,
\frac{(xJ_k^{\epsilon})_{\nu}}{x_{\nu}}\leq \lambda_{k-1}+\alpha\\
\textrm{If } |\kappa|= k,
\frac{(xJ_k^{\epsilon})_{\kappa}}{x_{\kappa}}=0.
\end{cases}
\end{equation*}
Hence $\lep\leq \lambda_{k-1}+\alpha$ for that choice of $\epsilon$.

\bigskip

Let us go back to the bound proved in Theorem \ref{Teps}.
We can simplify further this bound, getting rid of the eigenvector.
We obtain the following nicer, but weaker version:

\medskip
\begin{corollary} Let $C$ be a $\delta$-code such that $m-\delta^2\leq \lambda_{k-1}$.
Then,

\begin{equation}\label{ineg}
|C|\leq \frac{4\sum_{|\kappa|=k}d_{2\kappa}a_{\kappa}}{m-\lambda_k}.
\end{equation}
\end{corollary}
\medskip

\noindent\proof If $C$ satisfies $\delta^2>m-\lambda_{k-1}$, since 
$\lambda_{k-1}<\lep$ for all non-negative $\epsilon$ (from Theorem
\ref{Teps} (i)), 
the bound of Theorem \ref{Teps} (iv) applies to $C$.
We get, using Cauchy-Schwartz inequality, and $\lep<\lambda_{k}$:

\begin{align*}
|C| &\leq \frac{1}{m-\lep}\frac{\big(\sum_{|\kappa|=k}v_{\kappa}^{\epsilon}(\epsilon_{\kappa}+a_{\kappa})\big)^2}{\sum_{|\kappa|=k}d_{2\kappa}^{-1}\epsilon_{\kappa}{v_{\kappa}^{\epsilon}}^2}\\
&\leq \frac{1}{m-\lambda_k}\sum_{|\kappa|=k}d_{2\kappa}\frac{(\epsilon_{\kappa}+a_{\kappa})^2}{\epsilon_{\kappa}}.
\end{align*}

The function $z\to \frac{(z+a)^2}{z}$ is minimized over $]0,+\infty[$ when
$z=a$. We obtain, with $\epsilon_{\kappa}=a_{\kappa}$, the announced bound.

%\eproof

\section{Asymptotic behavior of the largest eigenvalue $\lambda_k$ of  $T_k$}\label{sec5}

In this section, we compute the limit taken by $\lambda_k$ when
the quotient $n/k$ tends to some fixed value (Theorem \ref{lambda}).
This result is needed to pass to the asymptotic in the inequality
(\ref{ineg}) for the size of a Grassmannian code.

We first need some very explicit formulas for the coefficients
of the symmetric matrix ${J'}_k$ 
associated to the endomorphism $T_k$, in the orthonormal
basis 
\\
$\{\sqrt{\dk}\Pk, |\kappa|\leq k\}$.  From now on we change our usual
convention: if not specified, $\kappa$ is a partition of degree
$s$.
The diagonal coefficients of ${J'}_k$ are the same as the
ones of $J_k$, while the other coefficients, denoted by
${A'}_s[\kappa,\mu]$,
satisfy

\begin{equation*}
{A'}_s[\kappa,\mu]={A}_s[\kappa,\mu]\sqrt{\frac{\dk}{d_{2\mu}}}.
\end{equation*}

To start with, we gather some known results on the polynomials $\Ck$.

\subsection{Review of some properties of the polynomials $\Ck$}

The coefficients  {\scriptsize $\bb{\mu}{\ka}$} and $\binom{\mu}{\ka}$ are defined
respectively by the following properties:
\begin{align}\label{ck}
\sigma \Ck&=\sum_{|\mu|=s+1} \bb{\mu}{\ka} \Cmu\\ \label{ck1}
\epsilon \Ck&=\sum_{|\nu|=s-1}\binom{\kappa}{\nu} \Cnu 
\end{align}
and have the following explicit expressions:
\begin{align}\label{bin}
\bb{\ki}{\ka}&=\prod_{\substack{j=1\\j\neq i}}^m
\frac{2\kappa_i-2\kappa_j+j-i+1}{2\kappa_i-2\kappa_j+j-i}\\\label{bin1}
\binom{\ki}{\ka}&=(\kappa_i+1+\frac{m-i}{2})\prod_{\substack{j=1\\j\neq i}}^m \frac{2\kappa_i-2\kappa_j+j-i+1}{2\kappa_i-2\kappa_j+j-i+2}
\end{align}
while any other values are equal to zero (see \cite[Lemma 7.5.7]{M}, \cite{M}, \cite[Th 14.1]{JC}, \cite{La1}).
The polynomials $\Ck$ are intimately related to the decomposition of
$\Glm$-modules (\cite[Theorem 5.2.9]{GW}):
\begin{equation*}
\mathcal{R}(\Glm/\Om)=\oplus_{\kappa} F_m^{2\kappa}.
\end{equation*}
For later use, we settle the notation: $\delta_{\kappa}:=\dim(F_m^{\kappa})$ and we recall the formula (\cite{FH}):
\begin{equation}\label{dimFm}
\delta_{\kappa}:=\dim(F_m^{\kappa})=\prod_{1\leq i<j\leq m}\frac{\kappa_i-\kappa_j+j-i}{j-i}.
\end{equation}

\subsection{Formulas for the coefficients of the matrix ${J'}_k$}

\medskip
\begin{proposition}\label{ppp1} The matrix $B_s$ has the following properties:
\begin{enumerate}
\item $B_s[\kappa,\kappa']=0$ for all $\kappa\ne \kappa'$.
\item If $m\leq n/2$, 
\begin{align*}
2B_s[\kappa,\kappa]=&\sum_{i\in u(\kappa)}\binom{\ki}{\kappa}\bb{\ki}{\ka}\frac{2{\kappa}_i+m+1-i}{2{\kappa}_i+n/2+1-i}\\
&-\sum_{i\in d(\kappa)}\binom{\kappa}{\kii}\bb{\ka}{\kii}\frac{2{\kappa}_i+m-1-i}{2{\kappa}_i+n/2-1-i}.
\end{align*}
\item If $m=n/2$, $B_s[\kappa,\kappa]=m/2$.
\end{enumerate}
\end{proposition}
\medskip

\noindent\proof We recall that the coefficients $\beta_{\kappa,\nu}$ are
defined by:

\begin{equation*}
\Pk=\beta_{\kappa}\Ck+\sum_{\nu \mid \kappa>\nu
}\beta_{\kappa,\nu}\Cnu.
\end{equation*}
Inverting these relations, we obtain coefficients
$\alpha_{\kappa,\nu}$ such that

\begin{equation*}
\Ck=\alpha_{\kappa}\Pk+\sum_{\nu \mid \kappa>\nu
}\alpha_{\kappa,\nu}\Pnu.
\end{equation*}

Taking into account the formula (\ref{ck}),
we obtain:

\begin{equation*}
B_s[\kappa,\kappa']={\beta_{\kappa}}\sum_{|\mu|=s+1}\bb{\mu}{\ka}
\alpha_{\mu,\kappa'}+\sum_{|\nu|=s-1}{\beta_{\kappa,\nu}}
{\bb{\ka'}{\nu}}\alpha_{\kappa'}.
\end{equation*}

We use the following obvious relations:
$\alpha_{\kappa}\beta_{\kappa}=1$ and
$\alpha_{\kappa}\beta_{\kappa,\kappa'}+\beta_{\kappa'}\alpha_{\kappa,\kappa'}=0$
to rewrite

\begin{align}\label{bk}
B_s[\kappa,\kappa']=\frac{\beta_{\kappa}}{\beta_{\kappa'}}\Big(-\sum_{|\mu|=s+1}&\bb{\mu}{\ka}
\frac{\beta_{\mu,\kappa'}}{\beta_{\mu}}\\
&+\sum_{|\nu|=s-1}
{\bb{\ka'}{\nu}}\alpha_{\kappa'}
\frac{\beta_{\kappa,\nu}}{\beta_{\kappa}}\Big).
\end{align}

Let us assume first that $\kappa\neq \kappa'$. Since
$\bb{\mu}{\ka}$ is non zero only if $\mu=\ki$ for some index
$i$, and also $\beta_{\mu,\kappa'}$ is non zero only if $\mu={\kkj}$
for some index $j$,
at most one term in the first summation may be non zero, and the same
argument holds for the second summation. We only have to consider the
case when $\kappa'$ satisfies: for some indexes $i\neq j$,
$\kappa'_i=\kappa_i+1$ and $\kappa'_j=\kappa_j-1$. The remaining terms
in the expression of $B_s[\kappa,\kappa']$ correspond to
$\mu=\ki={\kkj}$ and $\nu=\kjj={\kkii}$.

Moreover, the coefficients $\beta_{\mu,\kappa}$ are calculated in
\cite{JC}, and in particular we have:
\begin{equation}\label{beta}
\frac{\beta_{\kappa,\kjj}}{\beta_{\kappa}}=-\frac{1}{2}\binom{\kappa}{\kjj}\frac{2\kappa_j+m-1-j}{2\kappa_j+n/2-1-j}.
\end{equation}

Replacing in (\ref{bk}) we have
\begin{align*}
B_s[\kappa,\kappa']=\frac{\beta_{\kappa}}{\beta_{\kappa'}}
&\frac{2\kappa_j+m-1-j}{2(2\kappa_j+n/2-1-j)}.\\
&\Big(
\bb{\ki}{\ka}\binom{\ki}{\kappa'}
-\binom{\kappa}{\kjj}
\bb{\ka'}{\kjj}
\Big).
\end{align*}

Combining (\ref{ck}) and (\ref{ck1}) in the obvious relation:
\begin{equation*}
(\epsilon\sigma-\sigma\epsilon) \Ck=m\Ck
\end{equation*}
leads to:
\begin{equation}\label{bibinom}
\sum_{|\mu|=s+1}\bb{\mu}{\ka}\binom{\mu}{\kappa'}
-\sum_{|\nu|=s-1}\binom{\kappa}{\nu}
\bb{\ka'}{\nu}
=\begin{cases}
0 \ \textrm{ if } \kappa\neq \kappa'\\
m \ \textrm{ if } \kappa=\kappa'
\end{cases}
\end{equation}

From (\ref{bibinom}) we can conclude that $B_s[\kappa,\kappa']=0$ when
$\kappa\neq \kappa'$.

When
$\kappa=\kappa'$, replacing (\ref{beta}) in (\ref{bk}) leads to the
formula (ii). If moreover $n=m/2$, taking
account of
(\ref{bibinom}) we obtain  $B_s[\kappa,\kappa]=m/2$.

%\eproof

We now give explicit formulas for the coefficients of ${J'}_k$:

\medskip
\begin{proposition}\label{coeff}
With the following notations: 
\begin{displaymath}
q_i:=2\kappa_i-i+m
\end{displaymath}
\begin{displaymath}
N:=n-2m
\end{displaymath}
\begin{displaymath}
D(x)=\frac{x^2}{x^2-1} \quad\textrm{ and }
\end{displaymath}
\begin{displaymath}
\left\{\begin{array}{ll} C(x)=\frac{(x+1)(x+N)}{(2x+N)(2x+N+2)} & x\neq 0\\ C(0)=\frac{1}{N+2}&\end{array}\right.
\end{displaymath}

we have the expressions:
\begin{align*}
B_s[\kappa,&\kappa]=\\
&\frac{m}{2}-\frac{N}{4}\sum_{i\in u(\kappa)}\big(\prod_{\substack{j=1\\j\neq i}}^mD(q_i-q_j+1)\big)
\frac{q_i+2}{2q_i+N+2}\\
&+\frac{N}{4}\sum_{i\in d(\kappa)}\big(\prod_{\substack{j=1\\j\neq i}}^mD(q_i-q_j-1)\big)
\frac{q_i}{2q_i+N-2}.
\end{align*}
\begin{align*}
{A'}_s[\kappa,\ki]=\Big(\big(\prod_{j\ne
  i}&D(q_i-q_j+1)D(q_i+q_j+N+1)\big).\\
&C(q_i)C(q_i+1)\Big)^{1/2}.
\end{align*}

\end{proposition}
\medskip

\noindent\proof For the calculation of $B_s$, we replace in Proposition
\ref{ppp1} (ii) the formulas (\ref{bin}) and (\ref{bin1}), and take
account of Proposition
\ref{ppp1} (iii).

In order to calculate ${A'}_s[\kappa,\ki]$, we have already seen that:
\begin{equation*}
A_s[\kappa,\ki]=\bb{\ki}{\ka}
\left(\frac{\beta_{\ki}}{\beta_{\kappa}}\right)^{-1}.
\end{equation*}

We need a formula for
$\left(\frac{\beta_{\ki}}{\beta_{\kappa}}\right)^{-1}$. Expressions for the leading coefficients 
of the polynomials $\Ck$ and $\Pk$ can be found in \cite{Macdo} and \cite{V}.
Putting them together we find:
\begin{equation*}
\frac{\beta_{\ki}}{\beta_{\kappa}}=
\prod_{\substack{j=1\\j\neq i}}^m \left(\frac{q_i+q_j+N}{q_i+q_j+N+1}\right)\frac{(2q_i+N)(2q_i+N+2)}{(q_i+N)(q_i+N+1)}
\end{equation*}
where the last fraction must be understood as $(N+2)/(N+1)$ when $q_i=0$.
Joined with (\ref{ck}), we obtain
\begin{align*}
A_s[\ka,\ki]=\prod_{\substack{j=1\\j\neq i}}^m &\left(\frac{q_i-q_j+1}{q_i-q_j}\right)\left(\frac{q_i+q_j+N+1}{q_i+q_j+N}\right).\\
&\frac{(q_i+N)(q_i+N+1)}{(2q_i+N)(2q_i+N+2)}.
\end{align*}
Next we use (\cite{FH}):
\begin{align*}
\frac{d_{2\ki}}{d_{2\ka}}=\prod_{\substack{j=1\\j\neq i}}^m
&\left(\frac{q_i-q_j+2}{q_i-q_j}\right)
\left(\frac{q_i+q_j+N+2}{q_i+q_j+N}\right).\\
&\frac{(2q_i+N+4)(q_i+N)(q_i+N+1)}{(2q_i+N)(q_i+1)(q_i+2)}
\end{align*}
where the last fraction must be understood as $(N+4)(N+2)/2$ when $q_i=0$,
and we obtain the announced formula for $A'[\ka,\ki]$.

%\eproof

\subsection{The limit of $\lambda_k$}

Now $n$ varies with $k$ so we rather denote by $T_k^{(n)}$ the endomorphism 
defined previously and $\lambda_k^{(n)}$ its largest eigenvalue.

\medskip
\begin{theorem}\label{lambda}
If $n/2k\to \ell$, while $n\to +\infty $ and $k\to +\infty$,

$$\lim \lambda_k^{(n)}=4\frac{\ell+1/m}{(\ell+2/m)^2}.$$
\end{theorem}
\medskip

\noindent\proof We give careful proofs in the cases $m=1$ and $m=2$, and will
be more sketchy  in the general case. 
As it was noticed previously, when  $m=1$ the eigenvalues are the
zeros of the Jacobi polynomials;  their asymptotic is calculated in
\cite{KL}, 
exploiting the differential equation for the Jacobi polynomials
and  Sturm's method. Another approach, using chain sequences, is used
in \cite{IL}. However, none of these methods seem to generalize
easily to the several variable case.
Our argument  will only use the fact that the matrix ${J'}_k^{(n)}$
is non-negative. More precisely, we use the following:

\medskip
\begin{lemma}\label{lnn}\cite{G} 
Let $J$ be a non-negative symmetric matrix of size $N$, with largest eigenvalue
$\lambda$.
\begin{enumerate}
\item For all $x\in \R^N$ with $x_i>0$, $\lambda\leq
  \max_{i}\frac{(xJ)_i}{x_i}$.
\item For all $x\in \R^N$ $x\neq 0$, $\lambda\geq \frac{(xJ)\cdot x}{x\cdot x}$.
\end{enumerate}
\end{lemma}
\medskip

\noindent {\bf The case $m=1$.} We recover from Proposition \ref{coeff}
the formulas:

\begin{align*}
2b_s&=1-\frac{(n-2)(n-4)}{(4s+n)(4s+n-4)}\\
{a'}_s&=\left(
\frac{(2s+1)(2s+2)(2s+n-2)(2s+n-1)}{(4s+n-2)(4s+n)^2(4s+n+2)}\right)^{1/2}
\end{align*}

From these expressions we see that both sequences are increasing with
$s$. 
Moreover, we see easily that if $s\sim k$, $b_s\sim
2\frac{\ell+1}{(\ell+2)^2}$, and ${a'}_s\sim \frac{\ell+1}{(\ell+2)^2}$.
Applying Lemma \ref{lnn} (i) with $x_s=1$ for all $s$ 
leads to:

\begin{equation*}
\lambda_k^{(n)}\leq {a'}_{k-1}+b_k+{a'}_k
\end{equation*}
and the right hand side tends to $4(\ell+1)/(\ell+2)^2$ when $n/2k$ tends to $\ell$.

We lower bound $\lambda_k^{(n)}$ using Lemma \ref{lnn} (ii) and a
choice of $x$ proposed in \cite{IL}: let $x$ be defined by:

\begin{displaymath}
\left\{\begin{array}{llc}
x_s& = 0 & 1\leq s\leq t:=k-\lfloor \sqrt{k} \rfloor+1\\
x_s& = 1 & t+1\leq s\leq k+1
\end{array}\right.
\end{displaymath}
so that $x_s=1$ on the $\lfloor\sqrt{ k} \rfloor$ last coordinates. Then, 

\begin{align*}
\frac{(x{{J'}_k^{(n)}})\cdot x}{x\cdot x} &\geq \frac{\sum_{s=t+1}^{k-1} ({a'}_{s-1}+b_s+{a'}_s)}{k-t+1}\\
&\geq \left(\frac{k-t-1}{k-t+1}\right)({a'}_{t}+b_{t+1}+{a'}_{t+1})
\end{align*}

Again, the right hand side tends to $4(\ell+1)/(\ell+2)^2$, hence the result.

\medskip
\noindent{\bf The case $m=2$}. From Proposition \ref{coeff}, we have,
setting $s:=\kappa_1+\kappa_2$ and $v:=\kappa_1-\kappa_2$:

{\small
\begin{align*}
B_s[\kappa,\kappa]= 1+ &\frac{(n-6)(n-4)}{8(4s+2n-6)}
\left(\frac{(4s+2n-4)^2}{(4\kappa_1+n)(4\kappa_2+n-2)}\right.\\
&\left.-\frac{(4s+2n-8)^2}{(4\kappa_1+n-4)(4\kappa_2+n-6)}\right)
\end{align*}
\begin{align*}
{A'}_s[\kappa,&\kappa^{(1)}]=\left(\frac{(2v+2)^2}{(2v+2)^2-1}\right)^{1/2}\left(\frac{(2s+n-2)^2}{(2s+n-2)^2-1}\right)^{1/2}.\\
&\left(\frac{(2\kappa_1+2)(2\kappa_1+3)(2\kappa_1+n-3)(2\kappa_1+n-2)}{(4\kappa_1+n-2)(4\kappa_1+n)^2(4\kappa_1+n+2)}\right)^{1/2}
\end{align*}
\begin{align*}
{A'}_s[\kappa,&\kappa^{(2)}]=\left(\frac{(2v)^2}{(2v)^2-1}\right)^{1/2}
\left(\frac{(2s+n-2)^2}{(2s+n-2)^2-1}\right)^{1/2}.\\
&\left(\frac{(2\kappa_2+1)(2\kappa_2+2)(2\kappa_2+n-4)(2\kappa_2+n-3)}{(4\kappa_2+n-4)(4\kappa_2+n-2)^2(4\kappa_2+n)}\right)^{1/2}
\end{align*}
}

One can verify that these coefficients are increasing with $s$ when
$v$ stays constant. This is easy to see for $B_s$, not so obvious for
the two others because the second term is decreasing while the last
big quotient is increasing.

In order to obtain a lower bound for $\lambda_k^{(n)}$ from Lemma
\ref{lnn} (ii), we choose $x=(x_{\kappa})$ with:
$x_{\kappa}=0,1$. We fix a number $V<k-\sqrt{k}$.
Let ${\mathcal K}_{V,s}$ be the set of the $V$ partitions of
degree $s$ with smallest $v=\kappa_1-\kappa_2$. Hence ${\mathcal
  K}_{V,s}=\{\kappa\mid |\kappa|=s, \kappa_2\geq \lfloor
\frac{s}{2}\rfloor -V+1\}$.
We set $x_{\kappa}=1$ iff $\deg(\kappa)\geq t:=k-\lfloor\sqrt{k}\rfloor+1$, and 
$\kappa\in {\mathcal K}_{V,|\kappa|}$. We need to avoid in ${\mathcal
  K}_{V,s}$ some
partitions, namely the ones with $v=0$ and the ones with $v$
maximal (for those partitions, some terms are either missing or are
equal to zero in $(x{J'}_k^{(n)})_{\kappa}$). Let this new set be
denoted by ${\mathcal K'}_{V,s}$.
We have, when $x_{\kappa}=1$,  $\kappa \in {\mathcal K'}_{V,s}$,
$|\kappa|\neq t,k$,

\begin{align*}
(x{J'}_k^{(n)})_{\kappa}=B_s[\kappa,\kappa]&+\sum_{i=1}^2 {A'}_s[\kappa,\ki]\\
&+\sum_{i=1}^2 {A'}_{s-1}[\kii,\kappa]
\end{align*}

In the expressions of ${A'}_s[\kappa,\ki]$ we can minor the first term
by $1$ ($v\neq 0$), then minor each term by its minimal value in
the sequence $v=cte$ to which it belongs. As was  mentioned
before, this minimal value is obtained when the degree is minimal,
i.e. when $s= t$ or $s=t+1$. We do the same for
${A'}_{s-1}[\kii,\kappa]$ and for $B_s[\kappa,\kappa]$.
Then we must consider the behavior when $s$ is constant of
$B_s[\kappa,\kappa]$, of:

\begin{align*}
{A^1}_s[\kappa]:=
&\frac{(2s+n-2)}{((2s+n-2)^2-1)^{1/2}}\cdot\\
&\left(
\left(C(q_1)C(q_1+1)\right)^{1/2}+
\left(C(q_2)C(q_2+1)\right)^{1/2}\right)
\end{align*}

and of the analogous expression ${A^2}_{s-1}[\kappa]$ corresponding to the last term.
These expressions are increasing with $\kappa_2$. Let
${B_{s,V}}^{min}$, ${{A^i}_{s,V}}^{min}$
be their minimal values in ${\mathcal K'}_{V,s}$.
For simplicity, we assume that
$$\min({B_{t,V}}^{min},{B_{t+1,V}}^{min})=
{B_{t,V}}^{min},$$ and the same for $A^1$, $A^2$.

We obtain:

\begin{align*}
&\frac{(x{{J'}_k^{(n)}})\cdot x}{x\cdot x} \geq \\
&\frac{(k-t-1)(V-2)}{(k-t+1)V}({B_{t,V}}^{min}+{{A^1}_{t,V}}^{min}+{{A^2}_{t,V}}^{min}).
\end{align*}

Now we let $n/2k$ tend to $\ell$. Since ${B_{t,V}}^{min}$ is obtained
at a partition essentially equal to $[t/2-V/2,t/2+V/2$, and since
  $t\sim k$, we see that ${B_{t,V}}^{min}$ tends to
  $2(\ell+1/2)/(\ell+1)^2$. For the same reason, ${{A^1}_{t,V}}^{min}$ and
  ${{A^2}_{t,V}}^{min}$ tend to $(\ell+1/2)/(\ell+1)^2$ (the parameter
  $V$ is still fixed at this stage). So we obtain
\begin{equation*}
\liminf \lambda_k^{(n)} \geq (1-\frac{2}{V})\cdot 4\frac{\ell+1/2}{(\ell+1)^2}.
\end{equation*}

Now we let $V$ tend to $+\infty$ to obtain the appropriate lower
bound.

\medskip
The second and last step obtains an upper bound for $\lambda_k^{(n)}$
from Lemma \ref{lnn} (i) with an appropriate choice of $x$. The choice 
$x_{\kappa}=1$ for all $\kappa$ is not good enough here because
$D(2v+2)^{1/2}+D(2v)^{1/2}\neq 2$. We need some $x_{\kappa}$ that
modify properly these factors. We choose $x_{\kappa}:=(2v+1)^{1/2}$
where $v=\kappa_1-\kappa_2$. We have

\begin{align*}
\frac{(x{{J'}_k^{(n)}})_{\kappa}}{x_{\kappa}}=B_s[\kappa,\kappa]&+\sum_{i=1}^2 {A'}_s[\kappa,\ki]\frac{x_{\ki}}{x_{\kappa}}\\
&+\sum_{i=1}^2 {A'}_{s-1}[\kii,\kappa]\frac{x_{\kii}}{x_{\kappa}}
\end{align*}

Let:

{\small
\begin{align*}
{A^\sharp}[\kappa,&\kappa^{(1)}]=\left(\frac{2v+2}{2v+1}\right)\left(\frac{(2s+n-2)^2}{(2s+n-2)^2-1}\right)^{1/2}\cdot\\
&\left(\frac{(2\kappa_1+2)(2\kappa_1+3)(2\kappa_1+n-3)(2\kappa_1+n-2)}{(4\kappa_1+n-2)(4\kappa_1+n)^2(4\kappa_1+n+2)}\right)^{1/2}\\
{A^\sharp}[\kappa,&\kappa^{(2)}]=\left(\frac{2v}{2v+1}\right)\left(\frac{(2s+n-2)^2}{(2s+n-2)^2-1}\right)^{1/2}\cdot\\
&\left(\frac{(2\kappa_2+1)(2\kappa_2+2)(2\kappa_2+n-4)(2\kappa_2+n-3)}{(4\kappa_2+n-4)(4\kappa_2+n-2)^2(4\kappa_2+n)}\right)^{1/2}\\
{A^\flat}[\kappa_{(1)},&\kappa]=\left(\frac{2v}{2v+1}\right)\left(\frac{(2s+n-4)^2}{(2s+n-4)^2-1}\right)^{1/2}\cdot\\
&\left(\frac{(2\kappa_1)(2\kappa_1+1)(2\kappa_1+n-5)(2\kappa_1+n-4)}{(4\kappa_1+n-6)(4\kappa_1+n-4)^2(4\kappa_1+n-2)}\right)^{1/2}\\
{A^\flat}[\kappa_{(2)},&\kappa]=\left(\frac{2v+2}{2v+1}\right)\left(\frac{(2s+n-4)^2}{(2s+n-4)^2-1}\right)^{1/2}\cdot\\
&\left(\frac{(2\kappa_2-1)(2\kappa_2)(2\kappa_2+n-6)(2\kappa_2+n-5)}{(4\kappa_2+n-8)(4\kappa_2+n-6)^2(4\kappa_2+n-4)}\right)^{1/2}
\end{align*}
}

Since
$$D(2v+2)^{1/2}\frac{(2v+3)^{1/2}}{(2v+1)^{1/2}}=\frac{2v+2}{2v+1},$$
and
$$D(2v)^{1/2}\frac{(2v-1)^{1/2}}{(2v+1)^{1/2}}=\frac{2v}{2v+1},$$
we have:

\begin{align*}
\frac{(x{{J'}_k^{(n)}})_{\kappa}}{x_{\kappa}}=B_s[\kappa,\kappa]&+\sum_{i=1}^2 {A^\sharp}[\kappa,\ki]\\
&+\sum_{i=1}^2 {A^\flat}[\kii,\kappa].
\end{align*}

This expression is increasing with $s$ when $v$ is fixed. When $s$ is
fixed, $B_s[\kappa,\kappa]$, $\sum_{i=1}^2 {A^\sharp}[\kappa,\ki]
$ and $\sum_{i=1}^2 {A^\flat}[\kii,\kappa]$ are maximal at $\kappa=[s/2,s/2]$ (we extend the functions to
partitions with real parts here). We obtain, with $\rho_k=[k/2,k/2]$,

\begin{align*}
\max_{\kappa}&\frac{(x{{J'}_k^{(n)}})_{\kappa}}{x_{\kappa}}\leq \\
& B_s[\rho_k,\rho_k]+\sum_{i=1}^2 {A^\sharp}[\rho_k,{\rho_k}^{(i)}]
+\sum_{i=1}^2 {A^\flat}[{\rho_k}_{(i)},\rho_k].
\end{align*}

The computation of these values shows that the right hand side tends
to $4(\ell+1/2)/(\ell+1)^2$ when $n/2k\to \ell$.

\medskip
\noindent{\bf The general case $\bf m>2$} works the same. For the lower
bound,
we use ${\mathcal  K}_{V,s}:=\{\kappa\mid |\kappa|=s, \kappa_m\geq \lfloor
\frac{s}{m}\rfloor -V+1\}$. The cardinality of ${\mathcal  K}_{V,s}$
only depends on $s \mod m$. We should avoid some partitions in
${\mathcal  K}_{V,s}$, namely the ones with some parts equal and the
ones with $\kappa_m=\lfloor
\frac{s}{m}\rfloor -V+1$. Their
number is negligible compared to the cardinality of ${\mathcal  K}_{V,s}$.
Then, we proceed in the same way as for $m=2$.

\smallskip
The upper bound is obtained with
$x_{\kappa}=(\delta_{2\kappa})^{1/2}$.
We have 
\begin{equation*} 
\prod_{\substack{j=1\\j\neq i}}^mD(q_i-q_j+1)^{1/2}\left(\frac{\delta_{2\ki}}{\delta_{2\kappa}}\right)^{1/2}
=\prod_{\substack{j=1\\j\neq i}}^m\frac{q_i-q_j+1}{q_i-q_j}.
\end{equation*}
hence
\begin{align*}
\frac{(x{{J'}_k^{(n)}})_{\kappa}}{x_{\kappa}}=B_s[\kappa,\kappa]&+\sum_{i=1}^m {A^\sharp}[\kappa,\ki]\\
&+\sum_{i=1}^m {A^\flat}[\kii,\kappa],
\end{align*}
where
\begin{align*}
{A^\sharp}[\kappa,&\ki]=\prod_{\substack{j=1\\j\neq
    i}}^m\frac{q_i-q_j+1}{q_i-q_j}\cdot\\
&\Big(\big(\prod_{\substack{j=1\\j\neq i}}^mD(q_i+q_j+N+1)\big)
  C(q_i)C(q_i+1)\Big)^{1/2}
\end{align*}
and similarly
\begin{align*}
{A^\flat}[\kii,&\ka]=\prod_{\substack{j=1\\j\neq
    i}}^m\frac{q_i-q_j-1}{q_i-q_j}\cdot\\
&\Big(\big(\prod_{\substack{j=1\\j\neq i}}^mD(q_i+q_j+N-1)\big)
  C(q_i-2)C(q_i-1)\Big)^{1/2}.
\end{align*}

We have the nice identity:
\begin{equation*}
\sum_{i=1}^m \prod_{\substack{j=1\\j\neq i}}^m\frac{q_i-q_j+1}{q_i-q_j}=m.
\end{equation*}
We do not have a reference for this last identity, so we give an
argument here: from (\ref{dimFm}),
\begin{equation*}
\frac{\dim F_m^{\ki}}{\dim F_m^{\kappa}}=\prod_{\substack{j=1\\j\neq i}}^m\frac{\kappa_i-\kappa_j+j-i+1}{\kappa_i-\kappa_j+j-i}.
\end{equation*}
We obtain the demanded identity as the equality of the dimensions in
the following decomposition of $\Glm$-modules (Pieri's rule, \cite{GW}):
\begin{equation*}
F_m^{(1)}\otimes F_m^{\kappa}=\oplus_{i=1}^m F_m^{\ki}.
\end{equation*}

It turns out that the coefficients $B[\kappa,\kappa]$, $A^\sharp[\kappa,\ki]$
and $A^\flat[\kii,\ka]$ are increasing when $\kappa$ runs over
a sequence of the type $(\nu+ s[1,1,\dots,1])_{s\geq 0}$ (when $N$ is big enough), and that, on the space of partitions (with real parts) $\kappa$ of
fixed degree $k$, the maximum of the expressions
$B[\kappa,\kappa]$, $\sum_{i=1}^m A^\sharp[\kappa,\ki]$
and $\sum_{i=1}^m A^\flat[\kii,\ka]$ is attained at 
$\kappa=\rho_k=[k/m,k/m,\dots,k/m]$.

Moreover, it is easy to see that, when $n/2k\to \ell$,

\begin{displaymath}
\lim B_k[\rho_k,\rho_k]=2\frac{\ell+1/m}{(\ell+2/m)^2}
\end{displaymath}
\begin{align*}
\lim \left(\sum_{i=1}^m A^\sharp[\rho_k,\rho_k^{(i)}]\right)&=
\lim \left(\sum_{i=1}^m A^\flat[(\rho_k)_{(i)},\rho_k]\right)\\
&=\frac{\ell+1/m}{(\ell+2/m)^2}.
\end{align*}

%\eproof

\medskip
\begin{remark} One obvious consequence of Theorem \ref{lambda}
is that, for fixed $n$, the eigenvalue $\lambda_k^{(n)}$ runs over the 
whole interval $]0,m[$. Hence the bounds proved in section \ref{sec4} for the
  size of
grassmannian codes potentially cover all possible minimal distance.
\end{remark}
\medskip
\section{An asymptotic bound for the size of Grassmannian codes}\label{sec6}

We are now ready to take the limit when $n$ tends to $+\infty$ in the
inequality (\ref{ineg}), and prove Theorem \ref{Tass}.

\medskip
We are left with the estimate of
$\log(\sum_{|\kappa|=k}d_{2\kappa}a_{\kappa})/n$.

\medskip
\begin{lemma}
Let $\delta_{\kappa}:=\dim(F_n^{\kappa})$. If $n/2k\to \rho^{-1}\in \R$
while $n$ and $k$ tend to $+\infty$, 

\begin{align}\label{ineg2}
\limsup &\frac{1}{n}\log(\sum_{\substack{|\kappa|=2k\\\ell(\kappa)\leq m}}\delta_{\kappa}) \leq\\
&m\big((1+\rho)\log(1+\rho)-\rho\log(\rho)).
\end{align}
\end{lemma}
\medskip

\noindent\proof In the case $m=1$, $\delta_{2k}=\dim S_{2k}=\binom{n+2k-1}{2k}$  and it
is a classical result. The general case is probably well-known but
since we lack a reference, we give a proof here.
Let $\kappa$ be a partition of length at most $m$ and of degree $2k$,
that we extend to a partition with $n$ parts with an appropriate
number of zeros. From (\ref{dimFm}),
\begin{equation*}
\dim(F_n^{\kappa})=\prod_{1\leq i<j\leq n}\frac{\kappa_i-\kappa_j+j-i}{j-i}.
\end{equation*}
Since $\kappa_j=0$ when $j>m$, we have
\begin{equation*}
\dim(F_n^{\kappa})=\prod_{1\leq i<j\leq
  m}\frac{\kappa_i-\kappa_j+j-i}{j-i}
\prod_{1\leq i\leq m} \prod_{j>m}\frac{\kappa_i+j-i}{j-i}
\end{equation*}
We upper bound:
\begin{equation*}
\prod_{j>m}\frac{\kappa_i+j-i}{j-i}\leq \binom{n+\kappa_i-1}{\kappa_i}
\end{equation*}
and 
\begin{equation*}
\prod_{1\leq i<j\leq
  m}\frac{\kappa_i-\kappa_j+j-i}{j-i}\leq (2k+1)^{m^2}
\end{equation*}
to obtain 
\begin{align*}
\sum_{\substack{|\kappa|=2k\\\ell(\kappa)\leq m}}\delta_{\kappa} &\leq
(2k+1)^{m^2}\sum_{\substack{|\kappa|=2k\\\ell(\kappa)\leq m}}\left(\prod_{i=1}^m\binom{n+\kappa_i-1}{\kappa_i}\right)\\
&\leq (2k+1)^{m^2}\left(\sum_{s=0}^{2k} \binom{n+s-1}{s}\right)^m\\
&\leq (2k+1)^{m^2}\binom{n+2k}{2k}^m\\
\end{align*}
and we obtain the announced limiting result using the classical 
\begin{equation*}
\lim_{n/2k\to \rho^{-1}}
\frac{1}{n}\log\binom{n+2k}{2k}=(1+\rho)\log(1+\rho)-\rho\log(\rho).
\end{equation*}

%\eproof

From the three-term relation (\ref{3term}), specializing to $(1,\dots,1)$ we get
trivially\\
$a_{\kappa}\leq m$ and hence
$\sum_{|\kappa|=k}d_{2\kappa}a_{\kappa}\leq
m\sum_{\substack{|\kappa|=2k\\\ell(\kappa)\leq m}}\delta_{\kappa}$
(obviously $\dk\leq \delta_{2\kappa}$ since $V_n^{2\kappa}$ is
contained in $F_n^{2\kappa}$).

Then we only have to solve the equation, involving the limiting result
of Theorem \ref{lambda},
\begin{equation*}
s=\lim_{n/2k\to \rho^{-1}} \lambda_{k-1}^{(n)}=4\frac{\rho^{-1}+1/m}{(\rho^{-1}+2/m)^2}
\end{equation*}
which leads to
\begin{equation*}
\rho=\frac{m}{2}(-1+(1-\frac{s}{m})^{-1/2}).
\end{equation*}

%\eproof

\section{LP versus Hamming}\label{sec7}

In \cite{BN}, A. Barg and D. Nogin give an asymptotic bound for the size
of Grassmannian codes, derived from the so-called Hamming bound.
They prove, with the notations of Theorem~\ref{Tass}:

\medskip
\begin{theorem}\cite{BN}
\begin{equation}\label{BN}
\frac{1}{n}\log |C| \lesssim -m\log\left(\sqrt{1-\sqrt{\frac{s+m}{2m}}}\right)
\end{equation}
\end{theorem}
\medskip

It turns out that our bound (\ref{ass}) is better than (\ref{BN})
only when $s$ is small.\footnote{After this paper was submitted, the
  authors have further improved (\ref{BN}), see \cite{BN2}.}
The crossing point $s_0$ for the two bounds has the approximate value:

\begin{displaymath}
\begin{array}{c|c|c|c|c|c}
m&2&3&4&5&6\\
s_0& 1.4528 & 1.2714& 
1.1853&
1.1372&
1.1067\\
\end{array}
\end{displaymath}

\begin{displaymath}
\begin{array}{c|c|c|c}
7&8&9&10\\
1.0856&
1.0702&
1.0584&
1.0492\\
\end{array}
\end{displaymath}

\begin{figure}
\centering
\includegraphics{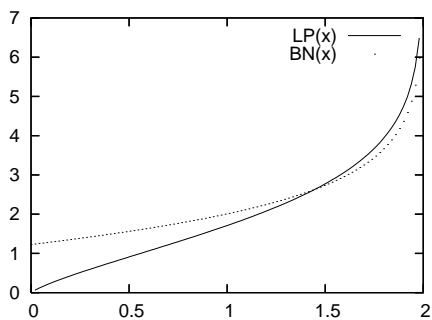}
\caption{LP and Hamming asymptotic bounds for $m=2$}\label{fig1}
\end{figure}

Figure \ref{fig1} plots the two bounds for $m=2$.

\section*{Aknowledgement}

We thank Pierre de la Harpe, Claude Pache,  Patrick Sol\'e, Gregory Kabatyanskiy 
and Vladimir Levenshtein for helpful discussions and suggestions. 
Part of this work
was done while the author was visiting Geneva University, supported by
the Swiss National Science Foundation.

\end{document}